\documentclass{article}
\pdfoutput=1

\usepackage{graphicx}
\usepackage{amsmath}
\usepackage{amssymb}
\usepackage{xfrac}
\usepackage{fix-cm}
\usepackage[affil-it]{authblk}

\author{David S. Abraham}
\affil{Department of Electrical and Computer Engineering, McGill University, Montr\'{e}al, QC H3A 0E9, Canada.}
\author{Alexandre Noll Marques}
\affil{Department of Aeronautics and Astronautics, Massachusetts Institute of Technology, Cambridge, MA 02139, USA.}
\author{Jean-Christophe Nave}
\affil{Department of Mathematics and Statistics, McGill University, Montr\'{e}al, QC~H3A~0B9, Canada.}

\begin{document}

\title{A Correction Function Method for the Wave Equation with Interface Jump Conditions}
\maketitle

\begin{abstract}
In this paper a novel method to solve the constant coefficient wave equation, subject to interface jump conditions, is presented. In general, such problems pose issues for standard finite difference solvers, as the inherent discontinuity in the solution results in erroneous derivative information wherever the stencils straddle the given interface. Here, however, the recently proposed Correction Function Method (CFM) is used, in which correction terms are computed from the interface conditions, and added to affected nodes to compensate for the discontinuity. In contrast to existing methods, these corrections are not simply defined at affected nodes, but rather generalized to a continuous function within a small region surrounding the interface. As a result, the correction function may be defined in terms of its own governing partial differential equation (PDE) which may be solved, in principle, to arbitrary order of accuracy. The resulting scheme is not only arbitrarily high order, but also robust, having already seen application to Poisson problems and the heat equation. By extending the CFM to this new class of PDEs, the treatment of wave interface discontinuities in homogeneous media becomes possible. This allows, for example, for the straightforward treatment of infinitesimal source terms and sharp boundaries, free of staircasing errors. Additionally, new modifications to the CFM are derived, allowing compatibility with explicit multi-step methods, such as Runge-Kutta (RK4), without a reduction in accuracy. These results are then verified through numerous numerical experiments in one and two spatial dimensions.\newline

\end{abstract}

\section{Introduction}
Wave problems with interface jump conditions are of primary importance in a wide variety of physical phenomena, forming the basis of fields as diverse as acoustics, elastodynamics, seismology, and electromagnetics. In general, it is found that wave discontinuities occur whenever two materials supporting different propagation velocities are in contact, or when singular source terms are present along a prescribed surface or curve, e.g surface charge distributions in electromagnetics. 


The goal of this paper is to develop a method that solves such wave problems with interface jump conditions to high order of accuracy. As discussed below, existing methods impose limits either in terms of accuracy (most are first or second order accurate), or in terms of the smoothness of the interface they can handle. The method presented here has no such limitations. However, as a first step in this development, attention is restricted to problems in which the domain is homogeneous, i.e. in which no material discontinuities occur. While material discontinuities play an important part in many physically significant problems, there are important problems that arise in the context of homogeneous domains, such as the case of surface charge distributions in electromagnetics. Furthermore, extensions of the current work to include material discontinuities are currently under investigation. These extensions explore the same approach adopted in \cite{[1]}, \cite{[14]}, and \cite{[28]}.

A brief overview of existing methods designed for interface wave problems is now presented. The Immersed Boundary Method (IBM), first proposed by Peskin \cite{[6],[27]} for the treatment of immersed boundaries in cardiac blood flow,  saw discontinuities re-expressed as singular source terms, which were then smeared over neighboring nodes, resulting in a method with first order accuracy. LeVeque and Li then proposed the Immersed Interface Method (IIM) \cite{[22],[25],[26]} for elliptic problems, in which the finite difference stencils are altered near the interface to compensate for discontinuities where they actually occur, yielding second order accuracy. The IIM was then further adapted to general hyperbolic equations by Zhang \cite{[3]}, and applied to specific physical problems in acoustics \cite{[4]}, as well as electromagnetics \cite{[2]}. Despite the successes of IIM, it nonetheless suffers from increased numerical dispersion near the interface, instabilities in the case of even modest differences in material characteristics \cite{[19]}, and is limited to second order accuracy.

Following these developments, Fedkiw et. al. introduced the Ghost Fluid Method (GFM) \cite{[20]} for multi-material flows. The GFM uses the concept of ``ghost'' nodes, in addition to standard nodes, within some small band surrounding the interface. The ghost variables therein serve as smooth extensions of the solution across the interface to the affected nodes, allowing for the compensation of discontinuities when calculating finite differences. These ghost corrections were then shown to be tantamount to an equivalent source term at the real nodes, allowing the use of existing solvers and preserving many desirable properties of the underlying schemes. The difficulty, however, lies in the accurate computation of these ghost terms at the affected nodes, and their effects on the order of accuracy. 

Building upon this concept of smooth solution extensions, Lombard and Piraud introduced a GFM inspired method known as the Explicit Simplified Interface Method (ESIM), specifically for both 1D and 2D acoustic and elastic waves \cite{[15],[19]}. The ESIM attempts to compute the smooth function extensions via a Taylor series expansion at the interface, coupled with the jump conditions, and in principle offers arbitrarily high accuracy. Unfortunately, however, it cannot treat the most general case of sharp interfaces, requiring curves be at least $C^1$.

In contrast, the present method is not only arbitrarily high order, but also capable of handling the most general interface geometries, including sharp curves of only $C^0$ smoothness. Based upon the recently developed Correction Function Method (CFM) for the Poisson equation \cite{[1]} (and which also saw application to the heat equation and to the Navier-Stokes equations \cite{[14]}), the method is itself inspired by the GFM. Rather than defining the required corrections at a discrete number of ghost nodes however, the CFM extends the concept to that of a correction function, similarly defined within some small region surrounding the interface. As a result, the correction function can be found to be governed by its own partial differential equation (PDE), which can, in principle, be solved to any order. In particular, by solving the defining PDE in a weak form via a least squares minimization procedure, the method is capable of solving a wide array of interface configurations and shapes, including those with problematic sharp points. 

As mentioned earlier, the present investigation focuses on the case of singular sources, in which it is assumed the material parameters are continuous within the problem domain (constant coefficient case). The resulting solution is therefore discontinuous across the interface defined as the surface (in 3D) or curve (in 2D) on which the sources reside. Despite these restrictions, it is emphasized that problems with continuously varying coefficients are equally permissible under certain circumstances, and that the generalization to material discontinuities is under development. Lastly, while the present treatment deals solely with static interfaces, in principle the generalization to moving interfaces is already possible with the current method. In this way, the CFM has the potential to allow discontinuous wave phenomena to be modeled in a comprehensive, accurate and versatile way.

\section{Theoretical Approach}
\subsection{Problem Definition}

\begin{figure}[hbtp]
\centering
\includegraphics[width=0.8\textwidth]{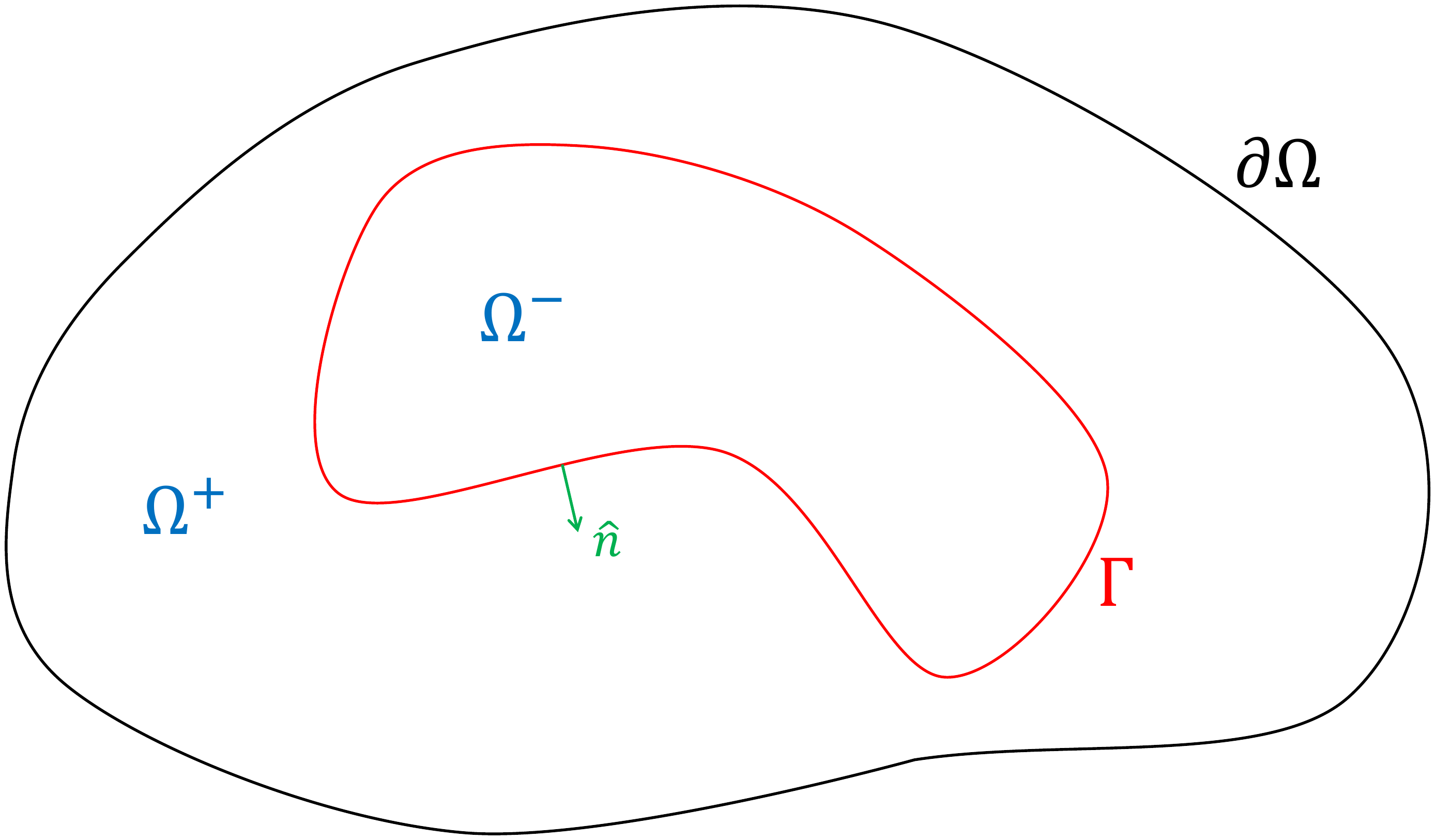}
\caption{Sample problem domain subdivided by an interface $\Gamma$.}
\label{fig1}
\end{figure}

The goal of this work is to solve the constant coefficient wave equation within a problem domain $\Omega$, in which the solution is discontinuous across an interface $\Gamma$ subdividing the domain into two regions, $\Omega^+$ and $\Omega^-$, as demonstrated in figure \ref{fig1}.

\newpage

For the entirety of this paper, $\vec{x} = (x_1, x_2, \ldots) \in \mathbb{R}^n$ is a spatial vector, $\nabla^2$ denotes the Laplacian defined as $\nabla^2 = \sum_{i=1}^n \frac{\partial^2}{\partial x_i^2}$, and lastly the normal derivative across the interface is given by $\frac{\partial u}{\partial n} = \hat{n}\cdot \nabla u$, in which $\hat{n}$ is an outward unit normal vector to the interface, as in figure \ref{fig1}.

The solution within region $\Omega^+$ shall be denoted as $u^+$, and likewise the solution in $\Omega^-$ as $u^-$. If $\vec{x} = \vec{x}_0 + \epsilon\hat{n}(\vec{x}_0)$, where $\vec{x}_0 \in \Gamma$, then a jump across the interface $\Gamma$ at $x_0$ is defined as:

\begin{equation}
[u(\vec{x}_0)]_\Gamma = \lim_{\epsilon \to 0^+} u(\vec{x}) - \lim_{\epsilon \to 0^-} u(\vec{x}).
\end{equation}

As such, the problem to be solved is as follows:

\begin{alignat}{4}
\nabla^2u(\vec{x},t)-\frac{1}{c^2}\frac{\partial^2u(\vec{x},t)}{\partial t^2} & = f(\vec{x},t) &\qquad& \text{in } \Omega\label{eqn1}\\
u(\vec{x},t) & = u(\vec{x}+\vec{L},t) &\qquad& \text{on } \partial\Omega\\
u(\vec{x},0) & = h(\vec{x}) &\qquad& \text{in } \Omega\\
\frac{\partial u}{\partial t}(\vec{x},0) & = k(\vec{x}) &\qquad& \text{in } \Omega\\
[u] & = \alpha(\vec{x},t) &\qquad& \text{on } \Gamma \label{eqn5}\\
\left[\frac{\partial u}{\partial n}\right] & = \beta(\vec{x},t) &\qquad& \text{on } \Gamma \label{eqn6}
\end{alignat}

\noindent in which $h(\vec{x})$, $k(\vec{x})$, and $f(\vec{x},t)$ are the usual prescribed initial conditions and forcing function, while $\alpha(\vec{x},t)$ and $\beta(\vec{x},t)$ represent the prescribed time dependent jumps in function value and normal derivative across $\Gamma$. In this paper the existence of periodic boundary conditions is assumed on the outermost domain boundary $\partial\Omega$, however in general the CFM can be applied regardless of the type of boundary condition used on $\partial\Omega$. 

\subsection{Motivating the Correction Function Method}
\begin{figure}[bhp]
\centering
\includegraphics[width=\textwidth]{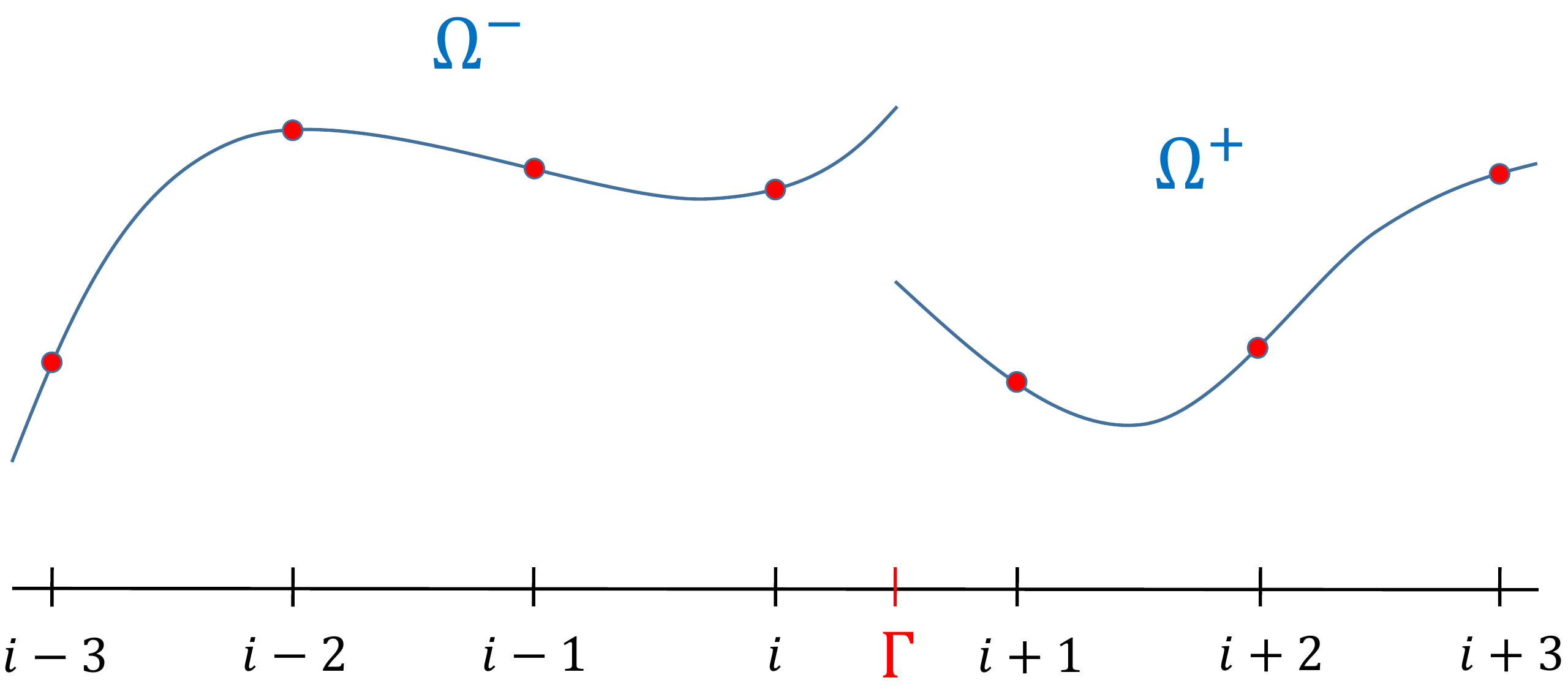}
\caption{Discretized one dimensional problem.}
\label{fig2}
\end{figure}

With the proposed problem defined, the use of the correction function is now motivated, demonstrating its basic principles, simplicity and benefits. As mentioned previously, solving the above type of problem with standard finite difference schemes yields erroneous results. To elucidate, consider a problem (without loss of generality) in one spatial dimension, as depicted in figure \ref{fig2}. Here it becomes obvious that any stencil straddling the interface makes use of the ``wrong'' function value on one side: $u^+$ rather than $u^-$ or vice versa. For example, suppose second derivative information at point $i$ was needed, equal to a forcing function $f(x)$, then:

\begin{equation}
\left.\frac{\partial^2 u}{\partial x^2}\right|_i \neq \frac{u^-_{i-1}-2u^-_i + u^+_{i+1}}{\Delta x^2} = f_i.
\end{equation}

\noindent This suggests searching for ways in which to ``correct'' the function on the other side of the interface, such that the above finite difference scheme yields the correct approximation. The Correction Function Method seeks to do exactly that.

\begin{figure}[tbhp]
\centering
\includegraphics[width=\textwidth]{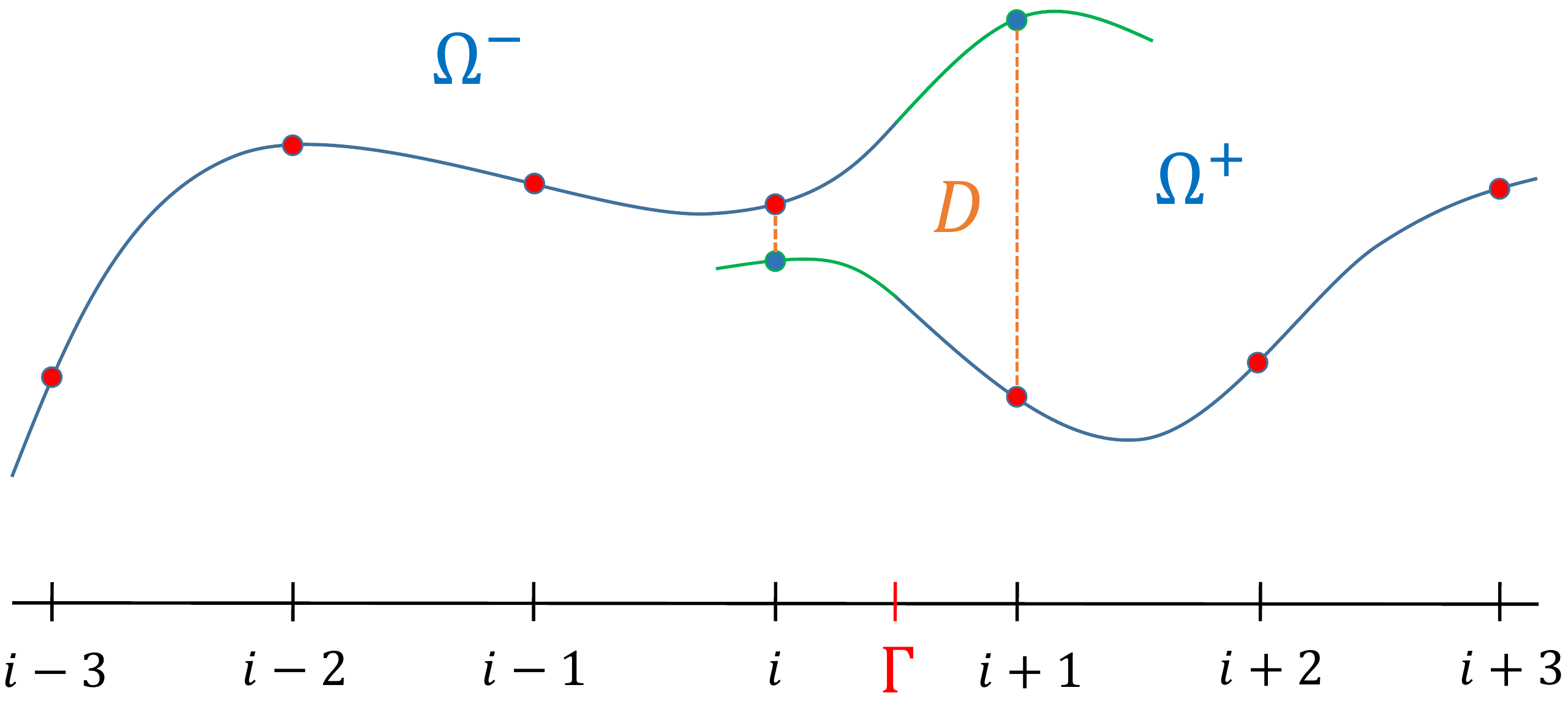}
\caption{Smooth extensions across the interface.}
\label{fig3}
\end{figure}

Assume for the moment that the solution on either side of the interface may be smoothly extended across the jump, as depicted in figure \ref{fig3}. From this figure, it is clear that the best points to use to obtain derivative information at nodes $i$ and $i+1$ belong to their respective extensions. To that effect, the correction function is defined as:

\begin{equation}
D = u^+(\vec{x},t) - u^-(\vec{x},t).
\end{equation}

\newpage

If it is equally assumed for the moment that the correction function $D$ so defined is known, then the stencil can be modified as follows:

\begin{equation}
\left.\frac{\partial^2 u}{\partial x^2}\right|_i \approx \frac{u^-_{i-1} - 2u^-_i + u^+_{i+1} - D_{i+1}}{\Delta x^2} = f_i. \label{eqn11}
\end{equation}

The inclusion of the correction function ensures that only points belonging to the continuous extension of the function are used, eliminating the discontinuity in the scheme. This method is very simple and, assuming the correction is known, is easy to implement, since the correction term can be included as an additional source term in existing solvers and methods:

\begin{equation}
\frac{u^-_{i-1} - 2u^-_i + u^+_{i+1}}{\Delta x^2} = f_i + \frac{D_{i+1}}{\Delta x^2}.\label{eqn12}
\end{equation}

\subsection{The Correction Function Method for the Wave Equation}
As was demonstrated in the previous section, if the correction function is known, solving the discontinuous interface problem becomes almost trivial. In other words, the burden has shifted from solving the original discontinuous problem, to finding an accurate approximation to the correction function $D$. In order to do so, a defining equation for $D$ is now derived, allowing for its solution wherever needed. To begin, the wave equation in each problem region is expressed separately:

\begin{equation}
\nabla^2u^+(\vec{x},t) - \frac{1}{c^2}\frac{\partial^2u^+(\vec{x},t)}{\partial t^2} = f^+(\vec{x},t)
\label{former}
\end{equation}
\begin{equation}
\nabla^2u^-(\vec{x},t) - \frac{1}{c^2}\frac{\partial^2u^-(\vec{x},t)}{\partial t^2} = f^-(\vec{x},t).
\label{latter}
\end{equation}

\noindent It is now assumed that $u^+$, $u^-$, $f^+$ and $f^-$ can each be smoothly extended across the interface, such that both sets of functions coexist in some small region around $\Gamma$. To ensure a leading error term $\mathcal{O}(\Delta x^4)$ in the solution $u$, it follows that the presumed extensions $f^+$ and $f^-$ need be at least $C^2$. With this assumption, equations (\ref{former}) and (\ref{latter}) may be subtracted within this area, yielding:

\begin{multline}
\nabla^2\left[u^+(\vec{x},t)-u^-(\vec{x},t)\right] - \frac{1}{c^2}\frac{\partial^2}{\partial t^2}\left[u^+(\vec{x},t)-u^-(\vec{x},t)\right] \\= f^+(\vec{x},t) - f^-(\vec{x},t) \equiv f_d(\vec{x},t).
\label{difference}
\end{multline}

\noindent In observing the above equation, the bracketed quantity is immediately recognized as the definition of the correction function. Making this substitution yields the desired defining equation for $D$:

\begin{equation}
\nabla^2 D(\vec{x},t) - \frac{1}{c^2}\frac{\partial^2D(\vec{x},t)}{\partial t^2} = f_d(\vec{x},t) \label{eqn16}.
\end{equation}

\noindent In like manner, the above procedure can be applied to the interface conditions \eqref{eqn5} and \eqref{eqn6} to obtain the needed interface conditions on $D$:

\begin{alignat}{4}
D(\vec{x},t) &= \alpha(\vec{x},t) &\qquad& \vec{x} \in \Gamma \label{eqn17}\\
\frac{\partial D(\vec{x},t)}{\partial n} &= \beta(\vec{x},t) &\qquad& \vec{x} \in \Gamma. \label{eqn18}
\end{alignat}

The above defining equations need only be solved within some small interval or band surrounding $\Gamma$, encapsulating the affected finite difference stencils. This is consistent with the derivation of \eqref{eqn16}, in which $D$ was only assumed to exist near the interface, where smooth extensions of $u^+$ and $u^-$ were also presumed to exist. 

It is worth noting that equations \eqref{former} - \eqref{difference} require that the same PDE be expressed in both regions. As a result, the method is currently inapplicable to problems in which the governing equations change across $\Gamma$, e.g. acoustic-elastic interfaces. Lastly, as mentioned, this paper not only assumes no material discontinuities are present, but also constant coefficients. However, the current method can, in theory, be used to solve problems with strictly continuously varying coefficients (again without material discontinuities). The procedure is a straightforward extension of the general CFM method, as detailed in \cite{[1]}, and will therefore not be addressed directly in this paper.

\subsection{Well-Posedness of the Wave CFM}
Having determined a defining partial differential equation for the correction function, the natural question of well-posedness arises, given that an arbitrary Cauchy problem is not always well-posed. Indeed, analysis of the elliptic Cauchy problem in \cite{[1]} demonstrated the ill-posedness of such a defining equation in a continuous setting. Though the ill-posedness of the Cauchy Poisson problem did not limit the applicability of the CFM, it was nonetheless an important consideration in the selection of an appropriate numerical scheme. This further underlines the need to understand the behaviour of equations \eqref{eqn16}, \eqref{eqn17} and \eqref{eqn18} in the present application.

To show that \eqref{eqn16} through \eqref{eqn18} are indeed well-posed, a method similar to that in \cite{[1]} is adopted. Assume, without loss of generality, that the interface is flat, and  that an orthogonal coordinate system $(\vec{y},d)$ is introduced on $\Gamma$, in which $\vec{y}$ spans the surface and $d$ measures the signed distance from the interface. Suppose now that a small perturbation to the original boundary conditions along the interface occurs. Taking the Fourier transform of the perturbation's $\vec{y}$ dependence, and solving the corresponding homogeneous wave equation, leads to a typical Fourier mode of the following form:

\begin{equation}
\phi_k = \mathrm{e}^{2\pi i(\vec{k}_y\cdot\vec{y}+k_d d\pm \omega t)}
\end{equation}

\noindent in which $k_d$ is a wavenumber along the $d$ axis, $k = |\vec{k}| = \sqrt{k_d^2 + \vec{k}_y\cdot\vec{k}_y}$ and $\omega = ck$. Here, in contrast to \cite{[1]} it is evident that such perturbations do not produce an exponential increase within the solution, only oscillation. As a result, the hyperbolic Cauchy problem defining the correction function is, in general, a well-posed problem with respect to the initial conditions. Given that the underlying equations are well behaved, numerical discretizations are therefore more likely to be well conditioned, yielding a much greater freedom in the selection of an appropriate numerical scheme.

\section{Implementation}
\subsection{Finite Difference and Time Marching Schemes}

Prior to discussing a solution strategy for the correction function itself, appropriate numerical schemes must be selected for the solution of the underlying wave equation. With both the underlying wave problem and the CFM Cauchy wave problem being well-posed, the freedom exists to choose any standard scheme. 

With this in mind, for the remainder of this paper, the non-compact, five point, fourth order accurate approximation for the spatial derivatives, defined as follows:

\begin{equation}
\left.\frac{\partial^2 u}{\partial x^2}\right|_i \approx \frac{1}{12\Delta x^2}(-u_{i-2} + 16u_{i-1} - 30u_i + 16u_{i+1} - u_{i+2}) \label{eqn24}
\end{equation}

\noindent will be used. Additionally, the fourth order accurate Runge-Kutta scheme (classical RK4) will be used to perform explicit time stepping. Each of these schemes was selected both for their simplicity, and ubiquitous nature.

\subsection{Computing the Correction Function}
\subsubsection{Defining the Correction Function Region}\label{sec3.2.1}
With the underlying schemes selected, the major task remains the calculation of the correction function via the defining equation. Despite having established that $D$ need only be computed within some small band around the interface which envelopes the finite difference stencil, a method to actually choose this region has not yet been devised, given that any number of regions may satisfy this criteria. For example, figure \ref{fig4} demonstrates a typical node in two dimensions whose stencils in $x$ and $y$ straddle the interface.

Clearly, the correction function will be required for some points within the region $\Omega^-$ to correct the stencil at the center point $(i,j)$. In this paper, only square regions will be considered as they greatly simplify the integrals required when $D$ is eventually solved in a least-squares sense (see section \ref{sec3.2.2}). Additionally, the amount of interface included in the region should be maximized, allowing for more information from $\Gamma$ to be included in the solution process, resulting in more accurate and consistent results. Lastly, while not as constrained as in the Poisson case, regions should not be made arbitrarily large due to accuracy and efficiency considerations. With these factors in mind, here a ``Node Centered Approach'' is adopted from \cite{[1]}, which may be summarized as follows:

\begin{figure}[b!]
\centering
\includegraphics[width=0.7\textwidth]{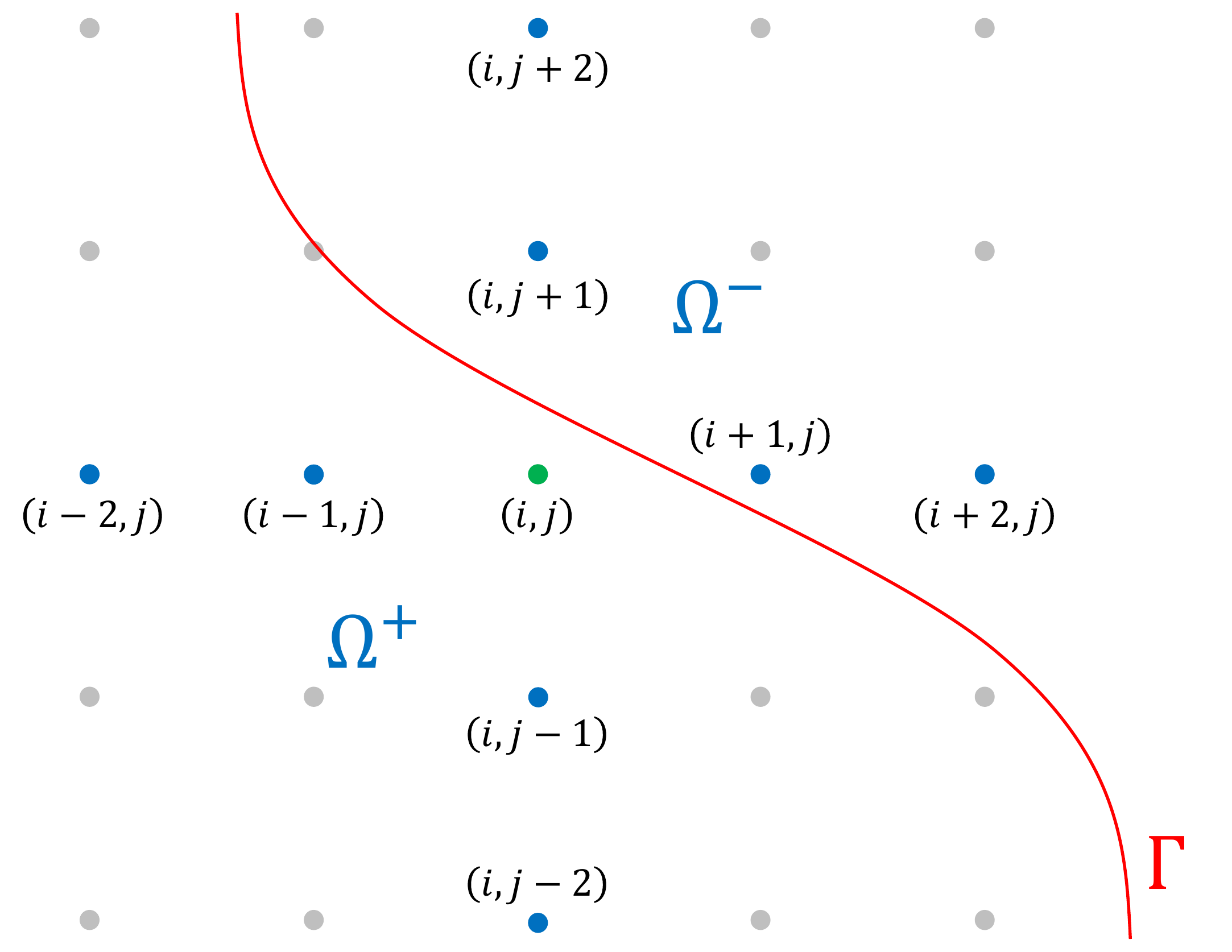}
\caption{Sample stencil in two dimensions straddling an interface.}
\label{fig4}
\end{figure}

\begin{enumerate}
\item Find the point $p_0$ along the interface which is closest to the node whose stencil is affected by $\Gamma$. To save time, $p_0$ may be roughly approximated.
\item Compute $\hat{t}_0$, the tangent vector to the interface at point $p_0$, as well as the normal vector $\hat{n}_0$ at point $p_0$. These two vectors now yield the directions of the principle diagonals of the square region.
\item The square region $\Omega_{i,j}$ is now defined as having sides of length $L$ and main diagonals aligned with $\hat{t}_0$ and $\hat{n}_0$.
\end{enumerate}

The value of $L$ in the above must be selected to ensure that all nodes associated with the center point $(i,j)$ of figure \ref{fig4} which lie on the opposite side of the interface ($\Omega^-$) are covered. This detail is crucial from an accuracy standpoint, as it ensures all corrected points associated with node $(i,j)$ share the same smooth error. If, for example, point $(i,j+1)$ was covered by one region and point $(i+2,j)$ by another, suboptimal accuracy would be observed, despite the correction functions in either region agreeing to leading order. The non-smooth nature of the leading error terms yield lower accuracy when computing the needed derivatives at $(i,j)$, similar to the loss of accuracy in calculating the gradients in \cite{[1]}. A good value for $L$ using stencil \eqref{eqn24} was found to be between $L = 3\sqrt{\Delta x^2 + \Delta y^2}$ and $L = 5\sqrt{\Delta x^2 + \Delta y^2}$, depending on the specific interface and grid geometry under consideration. Lastly, the results obtained in this paper assume the same value of $L$ for all regions, though a customized method in which each region selects a different, smallest, suitable value for $L$ may increase accuracy.

\begin{figure}[tbp]
\centering
\includegraphics[width=0.7\textwidth]{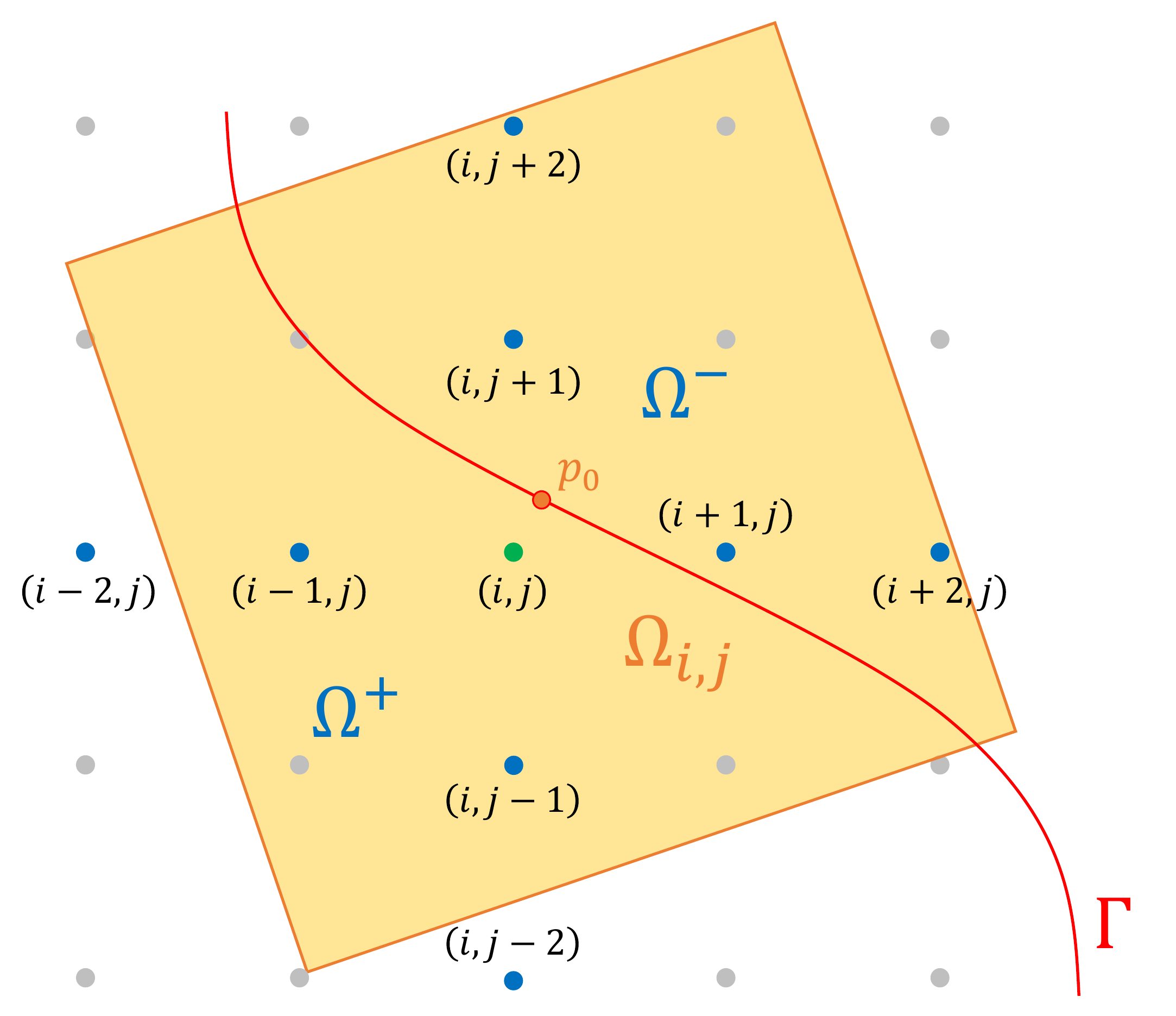}
\caption{Node centered approach to defining the correction function region.}
\label{fig5}
\end{figure}

\begin{figure}[tb]
\centering
\includegraphics[width=\textwidth]{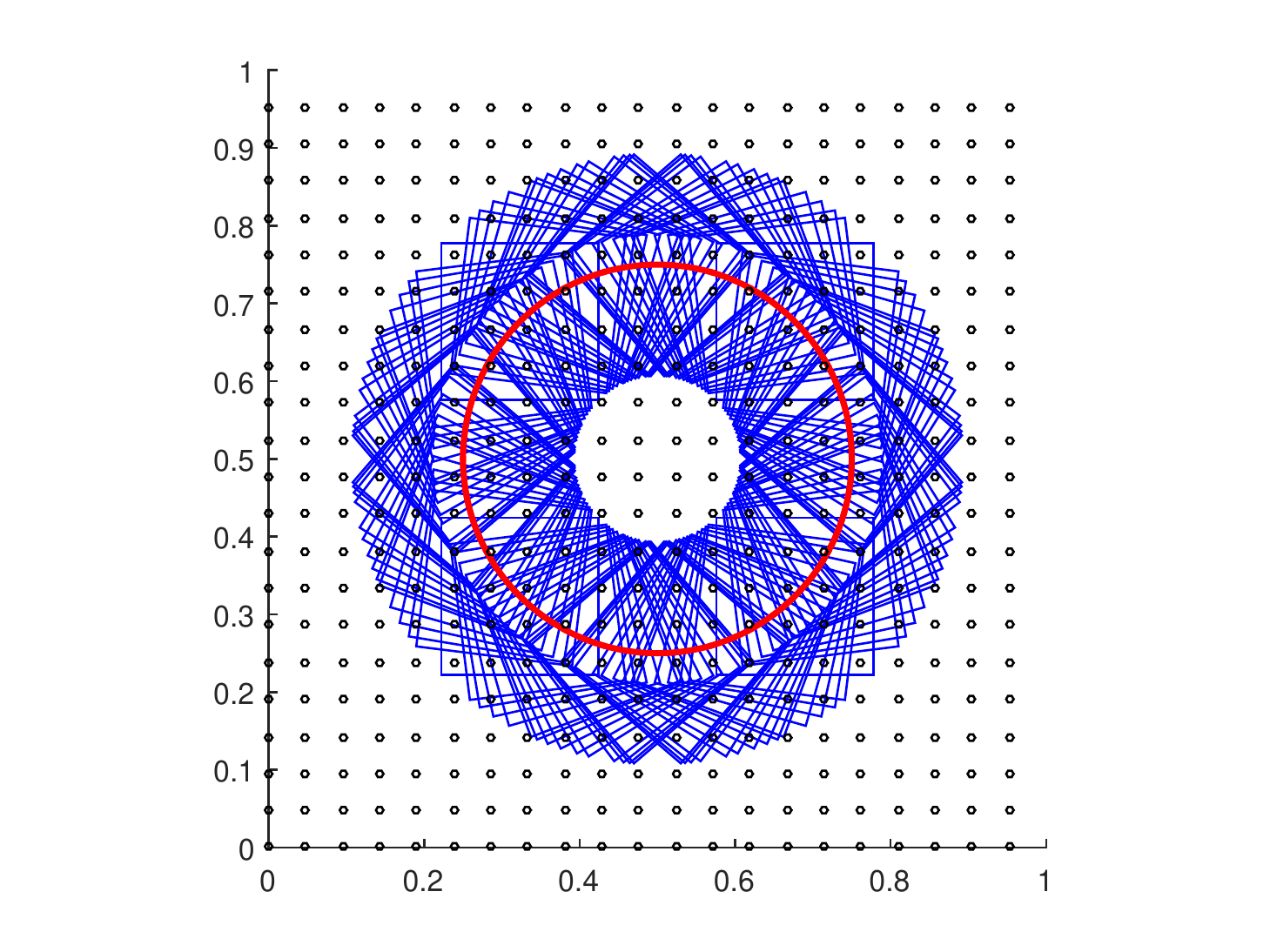}
\caption{Sample interface tiling for the calculation of the correction function in two dimensions.}
\label{fig6}
\end{figure}

By following this procedure for each region, the interface lies very close to the resulting main diagonals, ensuring that as much of it is encapsulated as is possible. This is exemplified in figure \ref{fig5} for the single node $(i,j)$. The process will then be repeated for each additional node whose stencil crosses the interface. Should the needed set of corrections for two different nodes both reside within a single region, a second region need not be defined, increasing efficiency, though this technique was not here employed. Each node will then make use of only those values of $D$ computed within its own region. A sample tiling is shown in figure \ref{fig6}, in which a circular interface is immersed in a 21 by 21 node mesh.

Lastly, since the wave problem has both spatial and temporal components, the required region should actually be a three dimensional space-time volume. However, given that all interfaces are here assumed static, it suffices to simply extend the squares generated by the above method into a rectangular prism of height $\Delta t$ (other multiples of $\Delta t$ may also be used). In this way, the correction function can be found within the given square spatial region for all times between $t_0$ and $t_0 + \Delta t$.

\subsubsection{Solving the Defining Equation} \label{sec3.2.2}
Given the vast number of possible interface/grid configurations which may arise, the technique of choice for solving the defining PDE for each sub-region is in a weak form, through the use of a least-squares minimization. This has the advantage of being robust and easily generalized. To do this, a suitable interpolant must first be selected with which $D$ may be approximated within each region.

Since a fourth order accurate finite difference stencil and time marching scheme are being employed, in order to preserve the leading fourth order accuracy of the overall scheme, the correction function must equally be known to at least fourth order. An initial observation of equation \eqref{eqn11} might lead to believe that $D$ is in fact required to sixth order, however this is not the case. Fourth order perturbations in $D$ can in effect be modeled as fourth order perturbations to the interface conditions encapsulated in equations \eqref{eqn5} and \eqref{eqn6}. These perturbations in turn will cause a slight fourth order alteration to the exact solution. Hence the correction function is only actually required to leading order, a postulate verified by numerical studies in section \ref{sec4}. As such, in this paper, tricubic Hermite interpolants have been selected in space and time, for problems in 2D.

With an appropriate temporal/spatial interpolant selected, the problem's weak form functional is defined as follows:

\begin{multline}
J_p = (l_c)^3\int_{t_0}^{t_0 + \Delta t}\!\!\int_{\Omega_{i,j}}\left[\nabla^2D - \frac{1}{c^2}\frac{\partial^2D}{\partial t^2} - f_d\right]^2\,\mathrm{d}\Omega_{i,j}\,\mathrm{d}t \\+ c_1\int_{t_0}^{t_0+\Delta t}\!\!\int_\Gamma [D-\alpha]^2\,\mathrm{d}\Gamma\,\mathrm{d}t + c_2(l_c)^2\int_{t_0}^{t_0+\Delta t}\!\!\int_\Gamma\left[\frac{\partial D}{\partial n}-\beta\right]^2\,\mathrm{d}\Gamma\,\mathrm{d} t.
\end{multline}

Clearly, the absolute minimum of this functional occurs for the exact solution of $D$, in which all integrands exactly vanish. The goal is to now replace $D$ within this functional by the tricubic Hermite interpolant representation, and minimize the resulting expression with respect to the expansion weights. The presence of scaling factors $l_c$ in the first and last terms ensures that all integrals scale similarly as the computational grid is refined, since each integrand has a different order of accuracy. The penalty coefficients $c_1$ and $c_2$ ensure that equal weighting is placed on all three conditions. In principle, they can be determined through analysis of the leading error term in $D$, but may be instead empirically determined through a coarse simulation, since they should not change significantly with refinement. In this paper, the integrals are carried out numerically by the use of Gaussian quadrature with six abscissae per dimension of integration. Combining the integration and minimization in this fashion in 2D yields a linear matrix system with 64 degrees of freedom, to be solved within every rectangular prism region, for each time interval. This results in the following global procedure for the current algorithm:

\begin{enumerate}
\item Create the tiling of regions $\Omega_{i,j}$ to cover the entirety of the interface, ensuring that each affected node has its required corrections contained within its region.
\item Solve the linear system for the correction function weights everywhere between the current step, $t_0$, and the subsequent step, $t_0+\Delta t$.
\item Advance the solution to $t_0 + \Delta t$ using the interpolated correction function in each region, in conjunction with the chosen finite difference schemes.
\item Set $t_0 = t_0 + \Delta t$ and repeat the process by returning to step 2, until the desired end time.
\end{enumerate}

\subsection{Accuracy and Stability}
\subsubsection{The Correction Function and RK4}
The classic RK4 scheme sees the solution marched forward in time by calculating the right-hand-side source term at several intermediate points, and taking a convex combination. Since RK4 is applied to a first order PDE, in the present case the second order wave equation must be recast into a coupled system of first order PDEs. Without loss of generality, this is expressed in one spatial dimension as follows:

\begin{align}
\frac{\partial u(\vec{x},t)}{\partial t} &= v(\vec{x},t)\label{eqn26}\\
\frac{\partial v(\vec{x},t)}{\partial t} &= c^2\left(\frac{\partial^2u(\vec{x},t)}{\partial x^2} - f(\vec{x},t)\right)\label{eqn27}
\end{align}

As such, the classic RK4 scheme takes on the following form for equation \eqref{eqn27}:

\begin{align}
k_{1,v} &= c^2\left(\frac{\partial^2 u}{\partial x^2} - f_n\right)\label{eqn28}\\
k_{2,v} &= c^2\left(\frac{\partial^2}{\partial x^2}\left\{u_n + \frac{\Delta t}{2}k_{1,u}\right\} - f_{n+\sfrac{1}{2}}\right)\label{eqn29}\\
k_{3,v} &= c^2\left(\frac{\partial^2}{\partial x^2}\left\{u_n + \frac{\Delta t}{2}k_{2,u}\right\} - f_{n+\sfrac{1}{2}}\right)\\
k_{4,v} &= c^2\left(\frac{\partial^2}{\partial x^2}\left\{u_n + \Delta t k_{3,u}\right\} - f_{n+1}\right)\label{eqn31}\\
v_{n+1} &= v_n + \frac{\Delta t}{6}(k_{1,v} + 2k_{2,v} + 2k_{3,v} + k_{4,v})
\end{align}

Clearly, from equations \eqref{eqn28} to \eqref{eqn31}, the finite difference stencil \eqref{eqn24} will need to be applied four separate times: once at $n$, twice at $n+\sfrac{1}{2}$, and once at $n+1$. Each of these applications will naturally require the use of the correction function for all stencils straddling the interface, at the given times. As such, a naive application of the CFM would see the interpolated correction function evaluated at the three time steps of interest, $D_n$, $D_{n+\sfrac{1}{2}}$ and $D_{n+1}$, for use in $k_1$ through $k_4$. Unfortunately, such a straightforward application results in sub-optimal accuracy on the order of $\Delta t^2$ rather than the expected $\Delta t^4$. To elucidate, consider equation \eqref{eqn29}, whereupon Taylor expanding the argument of the second derivative within $\Omega^+$ and $\Omega^-$ yields the following RK4 approximations to $u_{n+\sfrac{1}{2}}$:

\begin{align}
\hat{u}^-_{n+\sfrac{1}{2}} &\approx u^-_n + \frac{\Delta t}{2}k_1^- = u^-_n + \frac{\Delta t}{2}\left.\frac{\partial u^-}{\partial t}\right|_n\label{eqn33}\\
\hat{u}^+_{n+\sfrac{1}{2}} &\approx u^+_n + \frac{\Delta t}{2}k_1^+ = u^+_n + \frac{\Delta t}{2}\left.\frac{\partial u^+}{\partial t}\right|_n\label{eqn34}
\end{align}

In contrast, the correction function at this step is known to fourth order, and so it may equally be expressed as:

\begin{equation}
D_{n+\sfrac{1}{2}} \approx D_n + \frac{\Delta t}{2}\left.\frac{\partial D}{\partial t}\right|_n + \frac{\Delta t^2}{8}\left.\frac{\partial^2 D}{\partial t^2}\right|_n + \frac{\Delta t^3}{48}\left.\frac{\partial^3 D}{\partial t^3}\right|_n.\label{eqn35}
\end{equation}

Suppose now that derivative information is required at a node within $\Omega^+$ whose stencil crosses the interface (as in figure \ref{fig4}). The addition of $D$ would therefore be needed at points in $\Omega^-$ to generate $u^+$. Correspondingly, after adding equation \eqref{eqn35} to \eqref{eqn33} the following is obtained:

\begin{align}
\label{eqn36}
\begin{split}
\hat{u}^+_{n+\sfrac{1}{2}} &\approx (u^-_n + D_n) + \frac{\Delta t}{2}\left.\frac{\partial}{\partial t}(u^- + D)\right|_n + \frac{\Delta t^2}{8}\left.\frac{\partial^2 D}{\partial t^2}\right|_n + \frac{\Delta t^3}{48}\left.\frac{\partial^3 D}{\partial t^3}\right|_n\\
&\approx u^+_n + \frac{\Delta t}{2}\left.\frac{\partial u^+}{\partial t}\right|_n + \frac{\Delta t^2}{8}\left.\frac{\partial^2 D}{\partial t^2}\right|_n + \frac{\Delta t^3}{48}\left.\frac{\partial^3 D}{\partial t^3}\right|_n
\end{split}
\end{align}

The problem has now become evident: this last expression, equation \eqref{eqn36}, clearly differs from the required form given by equation \eqref{eqn34}. Consequently, the first two terms of equation \eqref{eqn36} will properly combine with their later counterparts in the final convex combination to yield fourth order accuracy, but the remaining terms involving $D$ will be carried through as extra error. In fact, if the CFM were applied in this way, each of the $k_2$ through $k_4$ steps would carry extra $D$ terms which would not combine to advance the solution ($k_1$ would not see additional error, since the value of $u$ used there is in fact known to fourth order, same as $D$). These error terms may be summarized as follows for each step of the process, in which only the leading order of the corrupting $D$ term has been kept:

\begin{align}
u^-_n + \frac{\Delta t}{2}k_1^- + D_{n+\sfrac{1}{2}} &= u^+_n + \frac{\Delta t}{2}k_1^+ + \mathcal{O}(\Delta t^2)\label{eqn37}\\
u^-_n + \frac{\Delta t}{2}k_2^- + D_{n+\sfrac{1}{2}} &= u^+_n + \frac{\Delta t}{2}k_2^+ + \mathcal{O}(\Delta t^3)\\
u^-_n + \Delta t k_3 + D_{n+1} &= u^+_n + \Delta t k_3^+ + \mathcal{O}(\Delta t^3)
\end{align}

These perturbations can in effect be modeled as slight alterations to the correction function, which as described previously is tantamount to changes in the interface conditions. With this point of view, the dominant second order error in equation \eqref{eqn37} will result in second order perturbations to the interface conditions, which in turn results in a second order discrepancy with the exact solution. This mismatch between the correction function and the RK4 intermediate steps has caused a reduction in the global truncation error from $\Delta t^4$ to $\Delta t^2$.

Nonetheless, an examination of equation \eqref{eqn36} suggests a way to compensate. By suppressing the second and third order terms in the Taylor expansions of $D$ within this expression, the discrepancy will no longer exist to leading order, and no sub-optimal error terms will accumulate. Through appropriate Taylor expansions of the spatial derivative arguments in equations \eqref{eqn28} through \eqref{eqn31}, the following altered forms of $D$'s expansion, consistent with RK4, may be derived:

\begin{align}
k_1 \mapsto \hat{D}_{n\phantom{+\sfrac{1}{2}}} &= D_n\\
k_2 \mapsto \hat{D}_{n+\sfrac{1}{2}} &\approx D_n + \frac{\Delta t}{2}\left.\frac{\partial D}{\partial t}\right|_n\\
k_3 \mapsto \hat{D}_{n+\sfrac{1}{2}} &\approx D_n + \frac{\Delta t}{2}\left.\frac{\partial D}{\partial t}\right|_n + \frac{\Delta t^2}{4}\left.\frac{\partial^2 D}{\partial t^2}\right|_n\\
k_4 \mapsto \hat{D}_{n+1\phantom{.5}} &\approx D_n + \Delta t\left.\frac{\partial D}{\partial t}\right|_n + \frac{\Delta t^2}{2}\left.\frac{\partial^2 D}{\partial t^2}\right|_n + \frac{\Delta t^3}{4}\left.\frac{\partial^3 D}{\partial t^3}\right|_n
\end{align}

To generate the temporal derivatives required in the above formulae, the tricubic interpolant may simply be differentiated and evaluated with the known weights at the current time step, $t_0$, and required nodes. In using this data in the appropriate way, the delicate combinations upon which RK4 are based are guaranteed to match, preserving accuracy.

While the above analysis focused solely on RK4, it is important to emphasize that a similar rationale and set of alterations may be applied to other time marching schemes as well, allowing them to also achieve their full potential. The CFM can therefore be used in a wide range of time marching schemes, as long as it is applied in a way consistent with the discretization being used.

\subsubsection{RK4 Stability for the Interface Wave Equation}\label{stab}
Having addressed accuracy in the previous section, attention is now turned toward the stability of the correction function augmented RK4 scheme. In equation \eqref{eqn12}, the correction function was reinterpreted as a source term on the right-hand-side of the standard finite difference scheme. Recasting the system of equations \eqref{eqn26} and \eqref{eqn27} into matrix form yields:

\begin{equation}
\frac{\partial}{\partial t} \begin{bmatrix} u\\v \end{bmatrix} = \begin{bmatrix} 0 & I\\ c^2 A & 0\end{bmatrix}\begin{bmatrix} u\\v\end{bmatrix} - \begin{bmatrix} 0\\c^2 f\end{bmatrix} \pm \begin{bmatrix} 0\\ \tilde{D}\end{bmatrix}
\end{equation}

\noindent where $A$ is the matrix representation of the finite difference scheme and $\tilde{D}$ is the equivalent source term derived from $D$, and is zero outside the interface region. Since the addition of the interface can simply be regarded as an additional source term, the stability of the scheme is independent of the interface. In other words, so long as the homogeneous (continuous) problem is stable, so too is the inhomogeneous interface problem. As such, given a bounded source term and interface conditions, the solution will remain bounded as well. The stability of the method presented here is therefore analyzed identically to that of standard RK4 and shares the exact same stability region. This observation is verified experimentally in the next section, and further demonstrates the robustness and versatility of the CFM in being incorporated into existing methods and solvers.

\section{Results}\label{sec4}
Numerical validation of the above derived method and analysis is now presented in both one and two spatial dimensions, for a variety of sample wave interface problems. In all cases, periodic boundary conditions have been applied to the problem domain's edges, with the interface(s) being described either by a single point (in one spatial dimension) or a curve. Additionally, in all cases the interfaces are assumed static with respect to time. Following these examples, the physically significant problem of electromagnetic radiation and shielding is demonstrated, as calculated in terms of the electric scalar potential.

\subsection{One Spatial Dimension}\label{line}
Rather than supplying all of the required conditions of equations \eqref{eqn1} through \eqref{eqn6}, the problem's exact solution $u$ and interface $\Gamma$ shall be specified, from which all the necessary conditions may be derived. For the one dimensional example, two interfaces have been specified, $\Gamma_1$ and $\Gamma_2$, with each being treated identically and independently over the domain $[0,1]$. Being points however, the node centered approach of section \ref{sec3.2.1} need not be used. Instead, a simple interval spanning two nodes on either side of the interface may be employed.

\noindent The exact solution is:

\begin{equation}
u(x,t) = 
\begin{cases}
\sin(2\pi x)\cos(2\pi t) &\text{if } 0 \leq x \leq \Gamma_1\\
2\sin(2\pi x)\cos(2\pi t) &\text{if } \Gamma_1 < x \leq \Gamma_2\\
\sin(2\pi x)\cos(2\pi t) &\text{if } \Gamma_2 < x \leq 1
\end{cases}
\end{equation}

\noindent where:

\begin{align}
\Gamma_1 &= \{x=0.3\}\\
\Gamma_2 &= \{x=0.7\}.
\end{align}

\begin{figure}[p]
\centering
\includegraphics[width=0.9\textwidth]{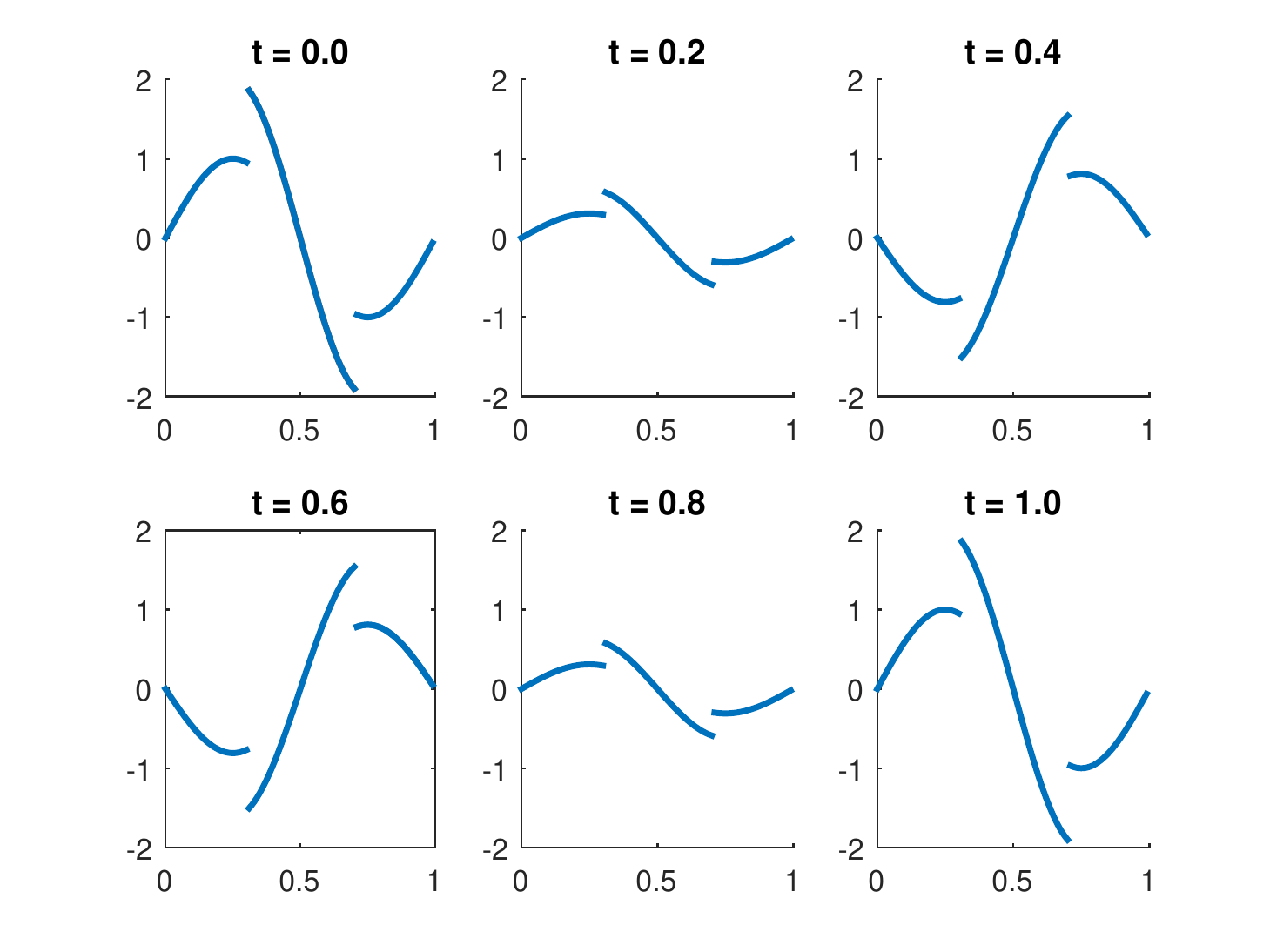}
\caption{Solution progression over time, demonstrating a crisp representation of the interface at all times.}
\label{fig7}
\end{figure}

The problem was solved for a full period, $t = 0$ to $t = 1$, maintaining the ratio $\Delta t = \Delta x$ at all times. Figure \ref{fig7} shows the solution at six intermediate time steps over one period, with $\Delta x = 0.01$. Importantly, the interface is seen to remain sharp and well defined, with no spurious oscillations developing.

\begin{figure}[p]
\centering
\includegraphics[width=0.9\textwidth]{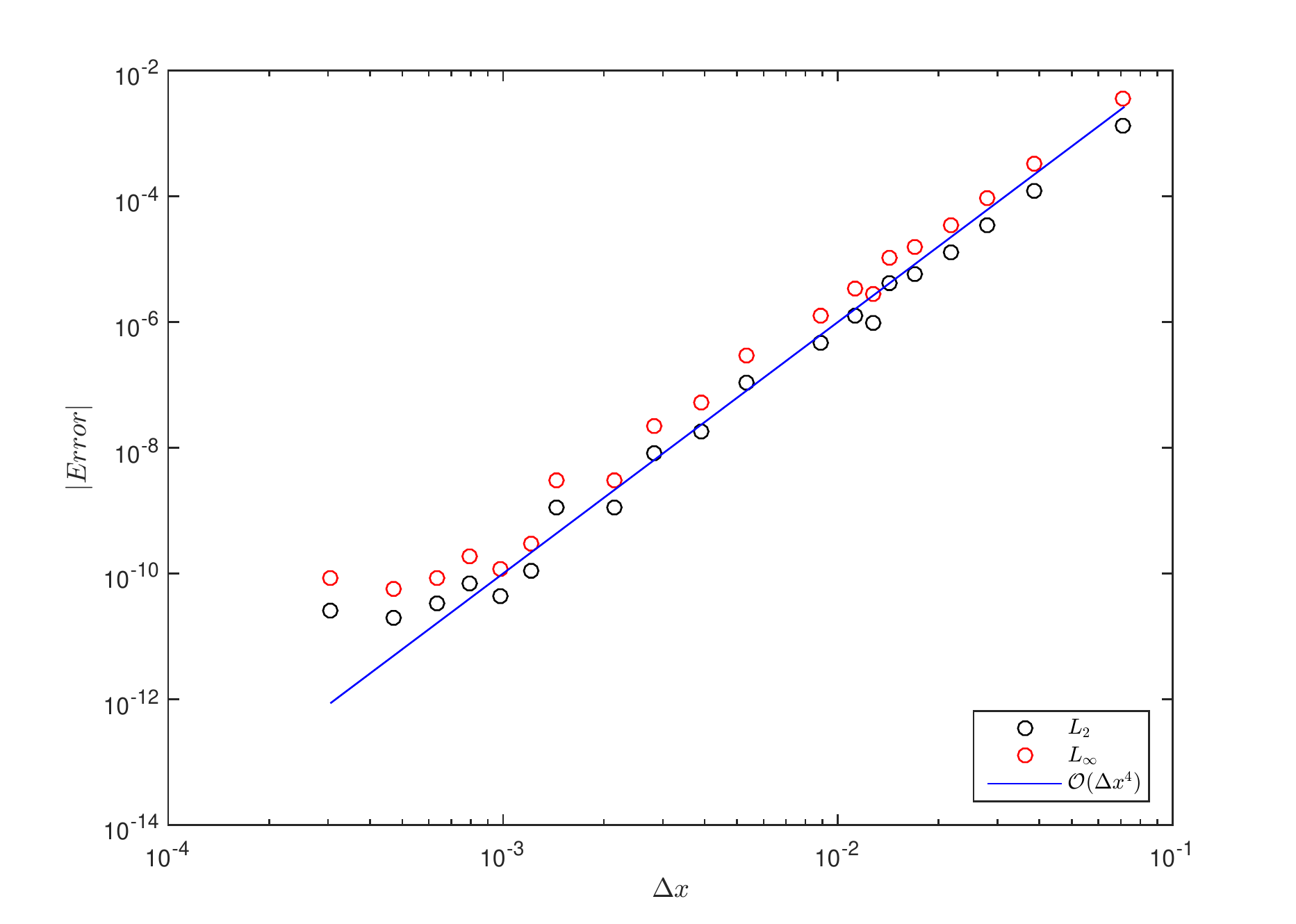}
\caption{1D Example -- Solution errors in $L_2$ and $L_\infty$ norms.}
\label{fig8}
\end{figure}

The error's convergence is evaluated in terms of the $L_2$ and $L_\infty$ norms, defined over all time steps and space as:

\begin{align}
\label{L2}
L_2 &= \frac{1}{\sqrt{N_x N_t}}\sqrt{\sum\limits_{x, t}(u-u_{ex})^2}\\
\label{Linf}
L_\infty &= \max_{x, t}|u - u_{ex}|
\end{align}

\noindent where $N_x$ is the number of nodes in the spatial direction, $N_t$ is the number of time steps, and $u_{ex}$ is the exact solution. Figure \ref{fig8} demonstrates the error as measured in the norms of equations (\ref{L2}) and (\ref{Linf}), as a function of the spatial node separation, $\Delta x$. The resulting trend demonstrates fourth order accuracy, as anticipated. Minor oscillations in the overall trend are to be expected, as the proportionality constant of the least squares error varies with the location of the interface within each region. As the grid is refined, the location of the interface shifts within each region, producing small jitters.

\subsection{Two Spatial Dimensions}
\subsubsection{Circular Interface}\label{circle}
The exact solution for the first 2D interface problem is defined over the region $[0,1]\times [0,1]$ and is as follows:

\begin{align}
\begin{split}
u^+(x,y,t) &= -2\sin(2\pi x)\sin(2\pi y)\cos(2\pi t)\\
u^-(x,y,t) &= \sin(2\pi x)\sin(2\pi y)\cos(2\pi t)
\end{split}
\end{align}

\noindent in which the interface is the circle centered at $(0.5,0.5)$ with radius $r = 0.25$. Additionally, given that the problem is now in two dimensions, the node centered approach of section \ref{sec3.2.1} was fully implemented.

\begin{figure}[p]
\includegraphics[width=\textwidth]{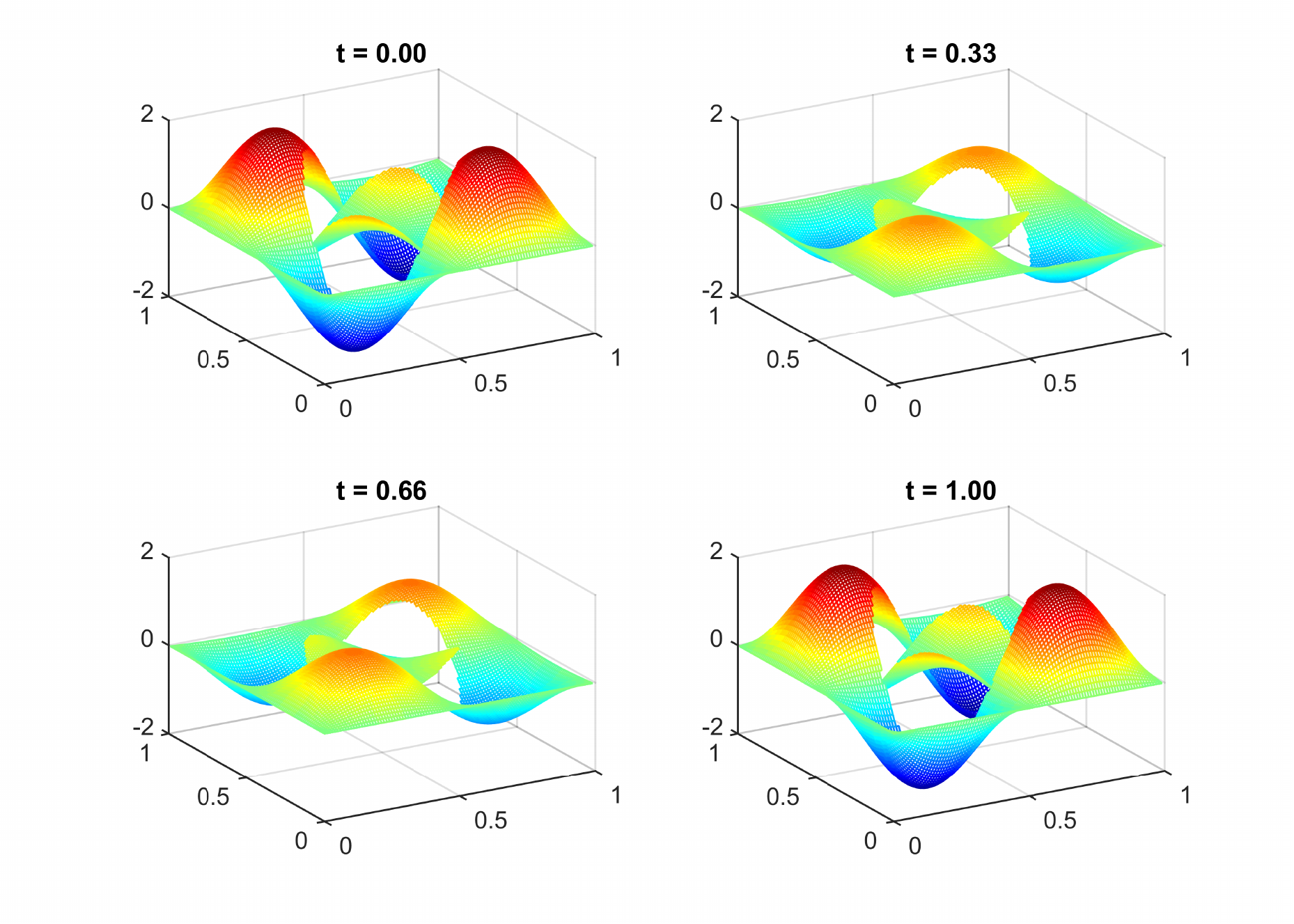}
\caption{Circular interface solution at various time steps.}
\label{fig9}
\end{figure}

\begin{figure}[p]
\includegraphics[width=\textwidth]{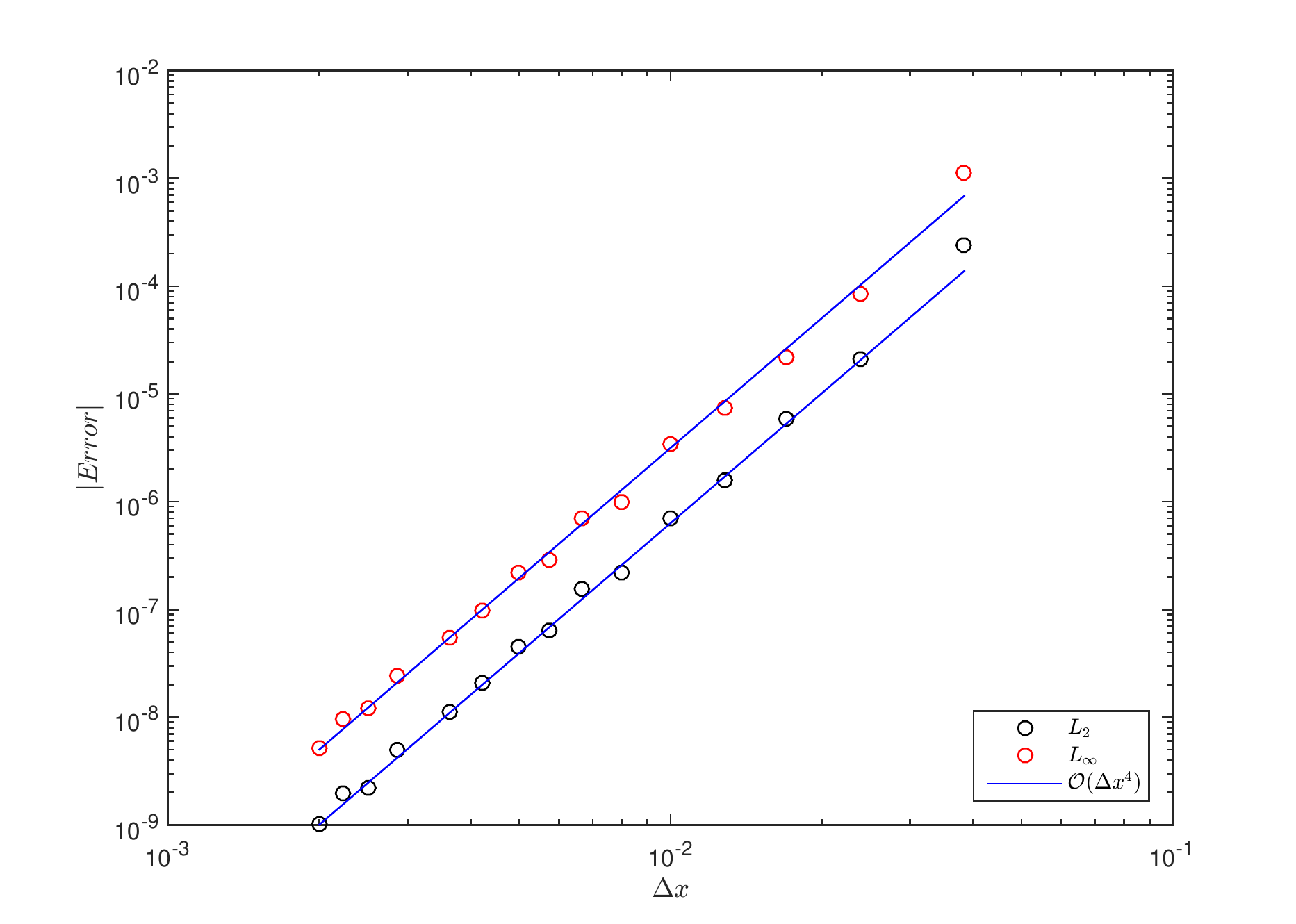}
\caption{Circular Interface -- Solution error in $L_2$ and $L_\infty$ norms.}
\label{fig10}
\end{figure}

The solution is plotted for four different times steps, over the interval $0 \leq t \leq 1$ with $\Delta x = \Delta t = 0.01$, in figure \ref{fig9}. Again the interface is seen to remain perfectly intact, with no perturbations. The error was computed for a range of $\Delta x$ values over one temporal period, and is plotted in figure \ref{fig10}. Again, fourth order accuracy can be observed in both the $L_2$ and $L_\infty$ norms.

\subsubsection{Star Interface}\label{star}
The second 2D example is similarly defined over the domain $[0,1]\times[0,1]$, with exact solution as follows:

\begin{align}
\begin{split}
u^+(x,y,t) &= 0\\
u^-(x,y,t) &= \mathrm{e}^{\pi x}\sin(3\pi y)\cos(2\pi t).
\end{split}
\end{align}

\begin{figure}[h]
\centering
\includegraphics[width=0.9\textwidth]{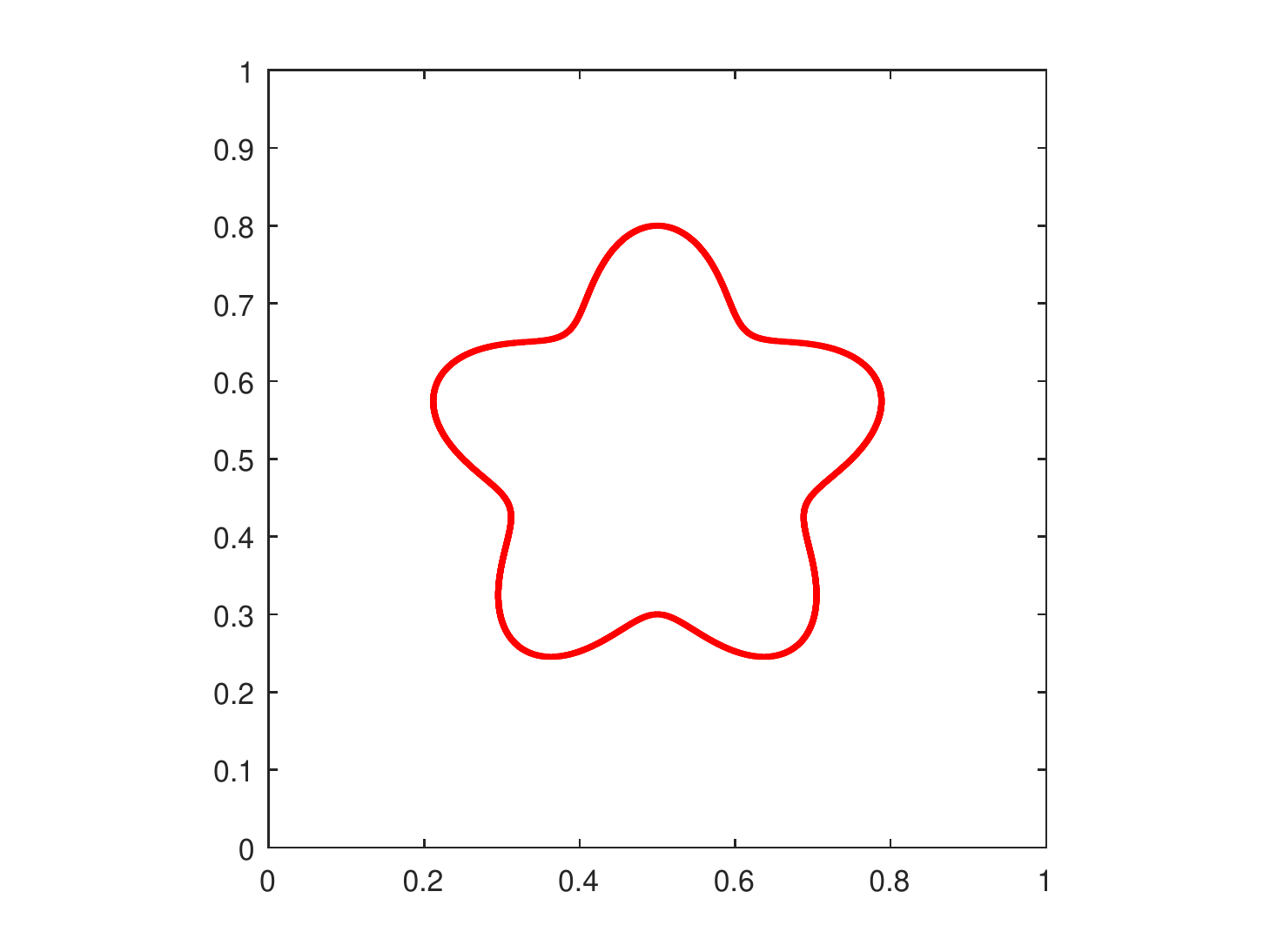}
\caption{Star-shaped interface in two dimensions.}
\label{fig11}
\end{figure}

Rather than being circular, the interface is now chosen to be a star shape, defined parametrically as:

\begin{align}
x &= (0.25 + 0.05\sin(5\theta))\cos(\theta) + 0.5\\
y &= (0.25 + 0.05\sin(5\theta))\sin(\theta) + 0.5
\end{align}

\noindent and shown in figure \ref{fig11}.

\begin{figure}
\includegraphics[width=\textwidth]{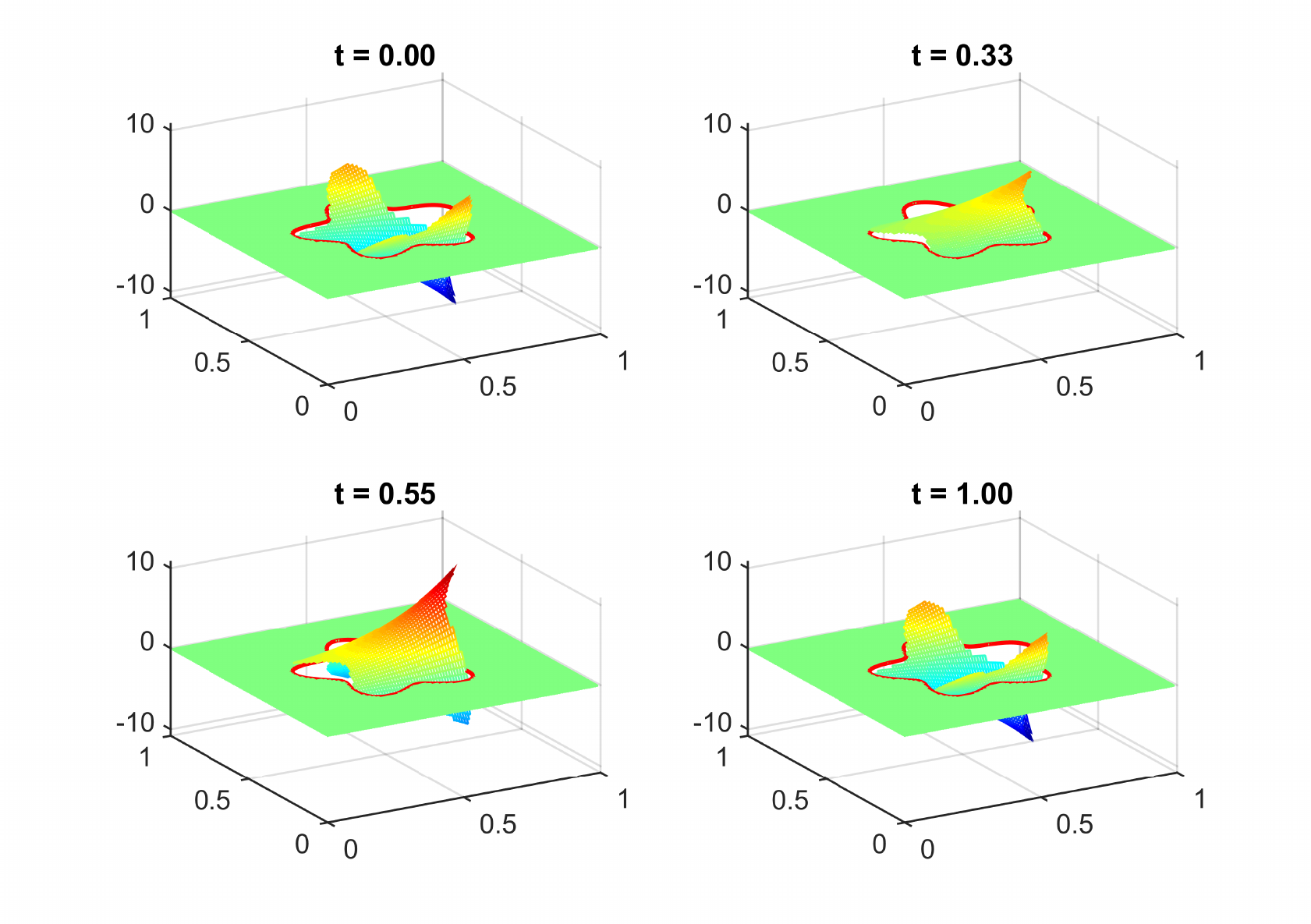}
\caption{Star interface solution at various time steps.}
\label{fig13}
\end{figure}

\begin{figure}
\includegraphics[width=\textwidth]{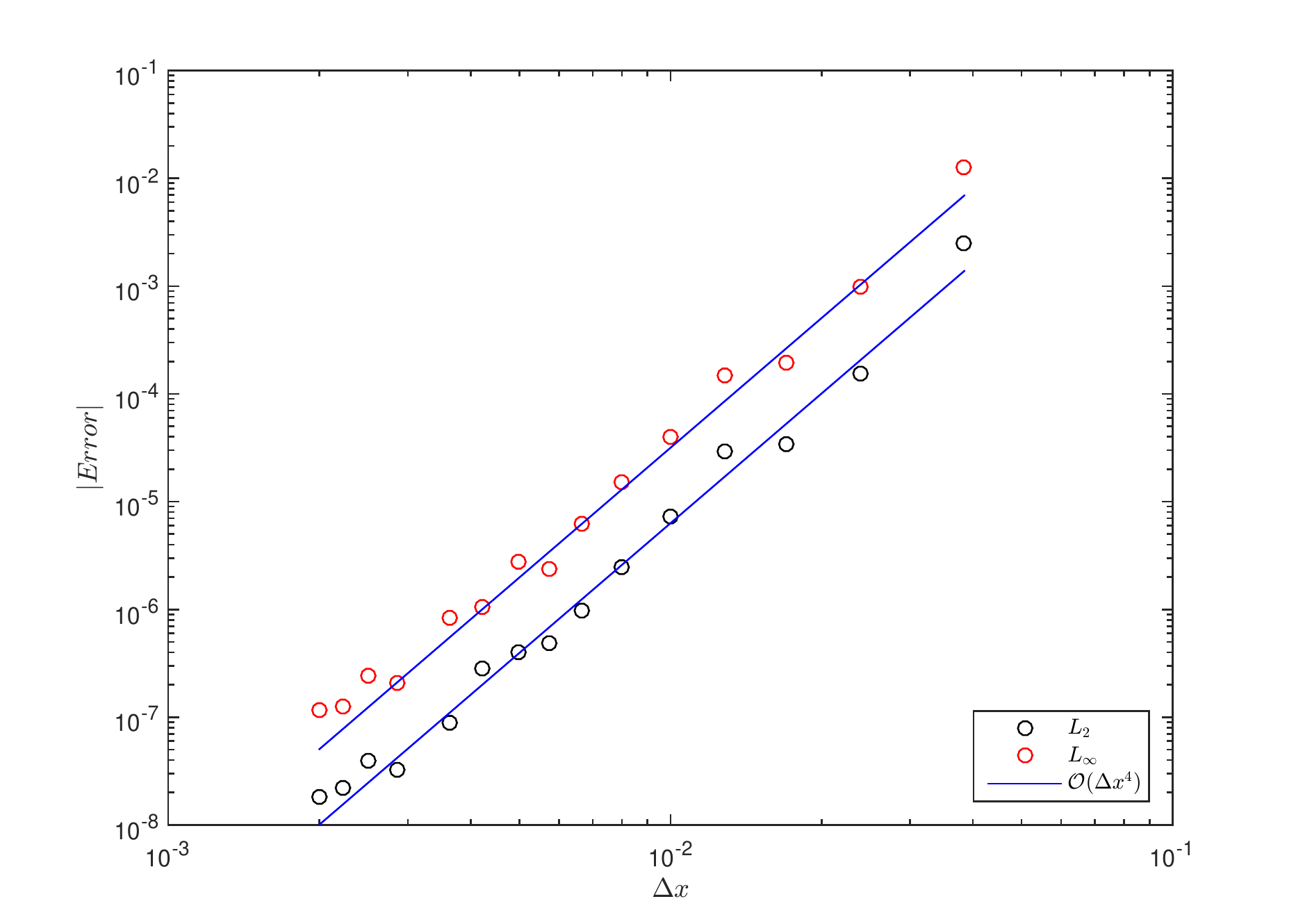}
\caption{Star Interface -- Solution error in $L_2$ and $L_\infty$ norms.}
\label{fig14}
\end{figure}

The solution is plotted in figure \ref{fig13} for several time steps over the interval $[0,1]$, with $\Delta x = \Delta t = 0.01$. Here, the sharpness of the discontinuity is particularly evident, given the zero solution outside the region (note the interface has been highlighted for clarity). Furthermore, the error convergence plot presented in figure \ref{fig14} demonstrates the expected fourth order accuracy in the $L_2$ and $L_\infty$ norms, over one temporal period.

\subsubsection{Sharp Interface}
The third 2D example is once again defined over the domain $[0,1]\times[0,1]$, and is determined by the solution:

\begin{align}
\begin{split}
u^+(x,y,t) &= 0.5\sin(2\pi x)\sin(2\pi y)\cos(2\pi t)\\
u^-(x,y,t) &= \mathrm{e}^{x + y}\cos(2\pi t).
\end{split}
\end{align}

\begin{figure}[h]
\centering
\includegraphics[width=0.9\textwidth]{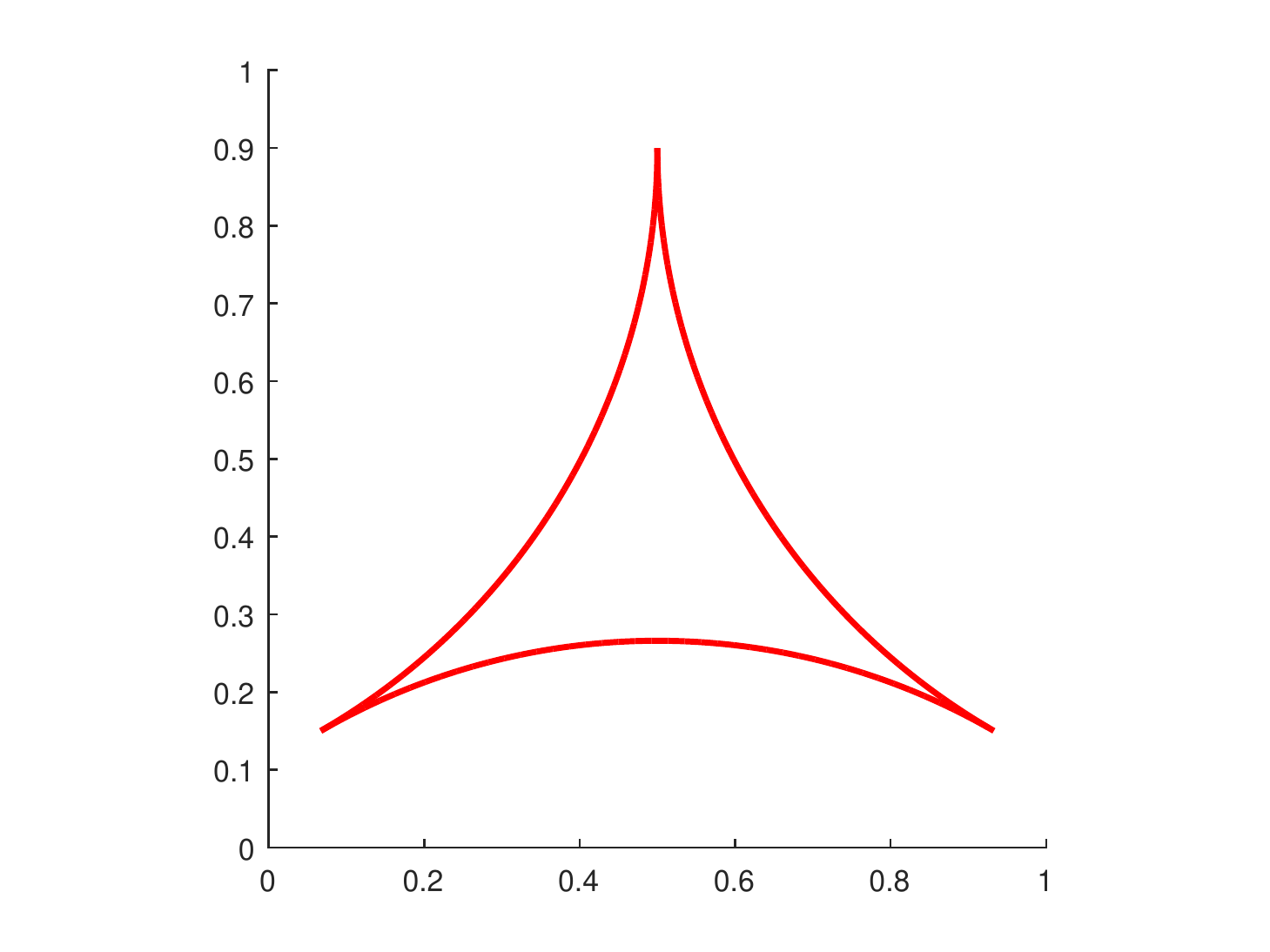}
\caption{Osculating circle interface in two dimensions.}
\label{fig12}
\end{figure}

The interface is now defined as bounding the region formed by three osculating circles, each with $r = \frac{\sqrt{3}}{2}$, and centered at $(0.5+\frac{\sqrt{3}}{2},0.9)$, $(0.5 - \frac{\sqrt{3}}{2}, 0.9)$ and $(0.5, -0.6)$. The interface is shown in figure \ref{fig12}.

\begin{figure}[p]
\includegraphics[width=\textwidth]{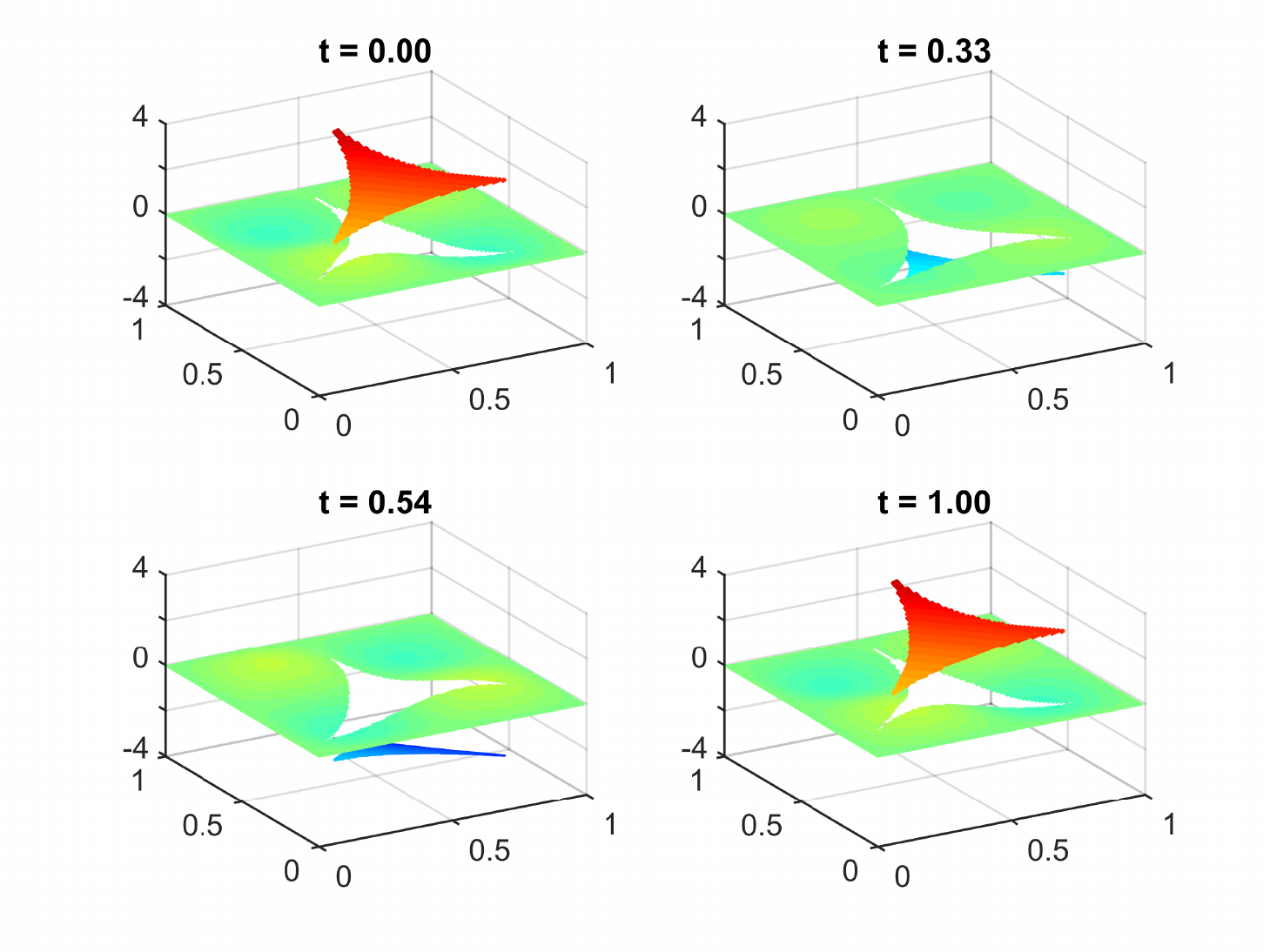}
\caption{Sharp interface solution at various time steps.}
\label{fig15}
\end{figure}

\begin{figure}[p]
\centering
\includegraphics[width=0.9\textwidth]{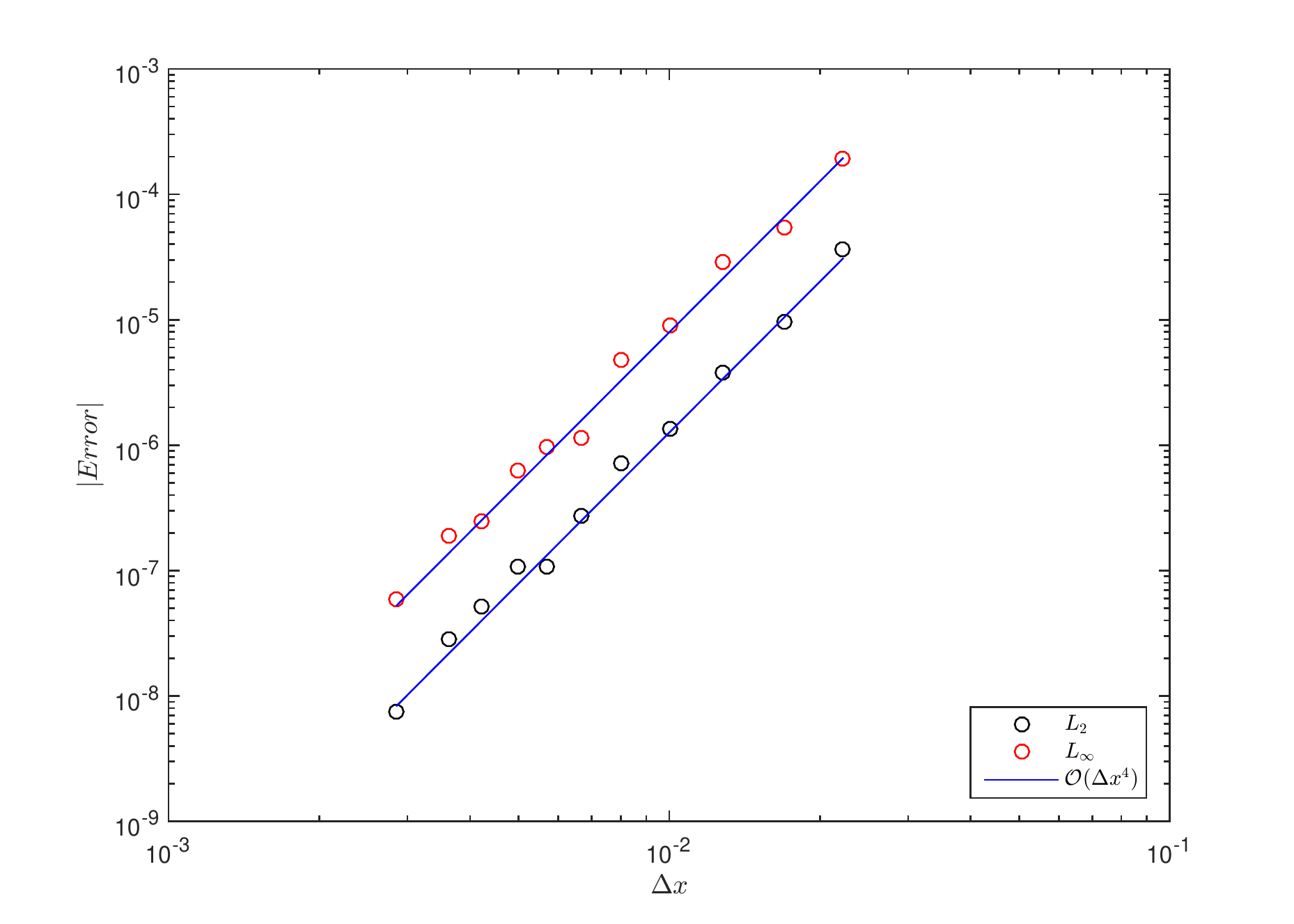}
\caption{Osculating Circles Interface -- Solution error in $L_2$ and $L_\infty$ norms.}
\label{fig16}
\end{figure}

Similarly to the previous two problems, the solution is shown for a variety of time steps in figure \ref{fig15}, again for $\Delta x = \Delta t = 0.01$. The error convergence plot is shown in figure \ref{fig16} for the $L_2$ and $L_\infty$ norms, computed over the time period $[0,1]$. Here, not only is the interface once again captured with high resolution, but the sharp points cause no added instability or reduction in accuracy, as evidenced by the 4th order trend of figure \ref{fig16}. Whereas other methods have traditionally struggled with such sharp interfaces, the CFM framework has handled them in a highly accurate and natural way.

\subsubsection{Maxwell's Equations of Electromagnetism}
Contrary to the previous three abstract examples, here a more physically meaningful problem (simplified slightly to two dimensions) is analyzed, in order to demonstrate the power and applicability of the developed method to real world problems.

Discontinuities are ubiquitous in the study of electromagnetic wave propagation and generation. Such waves often contain jumps when encountering dielectric interfaces, or when in the vicinity of charges and currents. In general, the behaviour of the electric ($\vec{E}$) and magnetic ($\vec{B}$) fields are governed by Maxwell's equations, a set of four vector partial differential equations. While many formulations solve for $\vec{E}$ and $\vec{B}$ directly, for the present method a more applicable approach is to recast the problem in terms of the vector and scalar potentials. A brief review of the pertinent electromagnetic theory is now presented, based upon \cite{[23]}.

The electric scalar potential $\phi$ and magnetic vector potential $\vec{A}$ are related to $\vec{E}$ and $\vec{B}$ by definition as follows:

\begin{align}
\vec{E} &= -\nabla\phi - \frac{\partial\vec{A}}{\partial t}\label{E}\\
\vec{B} &= \nabla\times\vec{A}\label{B}
\end{align}

\noindent from which it is possible to deduce a set of governing equations in terms of $\phi$ and $\vec{A}$ in free space:

\begin{align}
\nabla^2\phi + \frac{\partial}{\partial t}(\nabla\cdot\vec{A}) &= -\frac{\rho}{\epsilon_0}\label{phi}\\
\nabla^2\vec{A} - \frac{1}{c^2}\frac{\partial^2\vec{A}}{\partial t^2} -\nabla\left(\frac{1}{c^2}\frac{\partial\phi}{\partial t} + \nabla\cdot\vec{A}\right) &= -\mu_0 \vec{J}\label{A}
\end{align}

\noindent where $\rho$ is the volume electric charge density, $\vec{J}$ is the volume current density, $\epsilon_0$ and $\mu_0$ are the permittivity and permeability of free space, and $c$ is the speed of light in vacuum.

However, the divergence of $\vec{A}$ is not fixed by any of Maxwell's equations. This extra ``Gauge Invariance'' allows the divergence to be arbitrarily set. Two common choices are the Lorenz Gauge:

\begin{equation}
\nabla\cdot\vec{A} = -\frac{1}{c^2}\frac{\partial\phi}{\partial t}
\end{equation}

\noindent and the Coulomb Gauge:

\begin{equation}
\nabla\cdot\vec{A} = 0.
\end{equation}

In the Lorenz Gauge, the two equations \eqref{phi} and \eqref{A} become symmetric, each taking on the exact form of \eqref{eqn1}. Under the Coulomb Gauge, equation \eqref{phi} reduces to a Poisson problem for each time step, in which the method derived in \cite{[1]} may be used, with $\phi$ then generating an extra source term for $\vec{A}$. In either case, the following interface conditions can be derived for the constant coefficient problem \cite{[24]}:

\begin{align}
[\phi] &= 0 & \left[\frac{\partial\phi}{\partial n}\right] &= -\frac{\rho_s}{\epsilon_0}\label{jump2}\\
[\vec{A}] &= 0 & \left[\frac{\partial\vec{A}}{\partial n}\right] &= -\mu_0\vec{J}_s.
\end{align}

\noindent in which now $\rho_s$ and $\vec{J}_s$ represent charge and current densities confined to a surface.

These conditions, when coupled to the appropriate wave equation with consistent initial conditions, are ideally compatible with the method presented in this paper. While $\phi - \vec{A}$ formulations for Maxwell's equations in the time domain are not widespread within computational electromagnetics literature, they nonetheless have many attractive features. Methods which solve for the $\vec{E}$ and $\vec{B}$ fields directly in the time domain do have the advantage of not needing the post-processing of equations \eqref{E} and \eqref{B}. However, the potentials are in general much smoother functions, with the number of degrees of freedom reduced from six to four, and possibly less in the case of symmetries. Indeed, while not much literature is devoted to the topic, strong arguments for the consideration of potentials can be found in \cite{[8]}, and applications of $\phi$ and $\vec{A}$ to transmission line problems, Finite Element problems and eddy current calculations, can be found in \cite{[9]}, \cite{[13]} and \cite{[17]}, respectively.

With the above theoretical review complete, the calculation of the electric potential $\phi$ under the Lorenz Gauge is undertaken for an active shielding scenario. Given the complexity and general lack of closed form solutions for Maxwell's equations, here a manufactured solution is used to illustrate the general principle and also verify accuracy. The problem domain is defined as the square region $[-0.5,0.5]\times[-0.5,0.5]$, in which a $\phi$ plane wave is propagating diagonally. Centered at the origin, currents and charges are assumed to exist in the star-shaped pattern of section \ref{star}, actively canceling out the plane wave everywhere inside the region. As a result, the exact solution is as follows:

\begin{align}
\label{shield}
\begin{split}
\phi^+(x,y,t) &= \sin(2\pi[x+y]-\omega t)\\
\phi^-(x,y,t) &= 0
\end{split}
\end{align}

\noindent in which $\omega = 2\sqrt{2}\pi c$ rad/s, and $c$ is the speed of light in vacuum, 299792458~m/s.

Under these conditions, the surface charge and current distributions are precisely those required to shield the interior with respect to $\phi$, and may be explicitly computed through \eqref{shield} and \eqref{jump2}. It is worth noting that the solution presented in \eqref{shield} contains a discontinuity in the potential, which is not technically permissible given \eqref{jump2}. However, while \eqref{jump2} is the most common form in real-world applications, it is theoretically possible for $\phi$ to be discontinuous if more exotic charge configurations, such as an infinitely thin dipole layer \cite{[23]}, are present. While such configurations are rare, they are nonetheless assumed to exist in this case, for simplicity. A more practical implementation without the goal of verifying accuracy and performance would have no trouble applying equation \eqref{jump2} as stated.

\begin{figure}[p]
\centering
\includegraphics[width=0.9\textwidth]{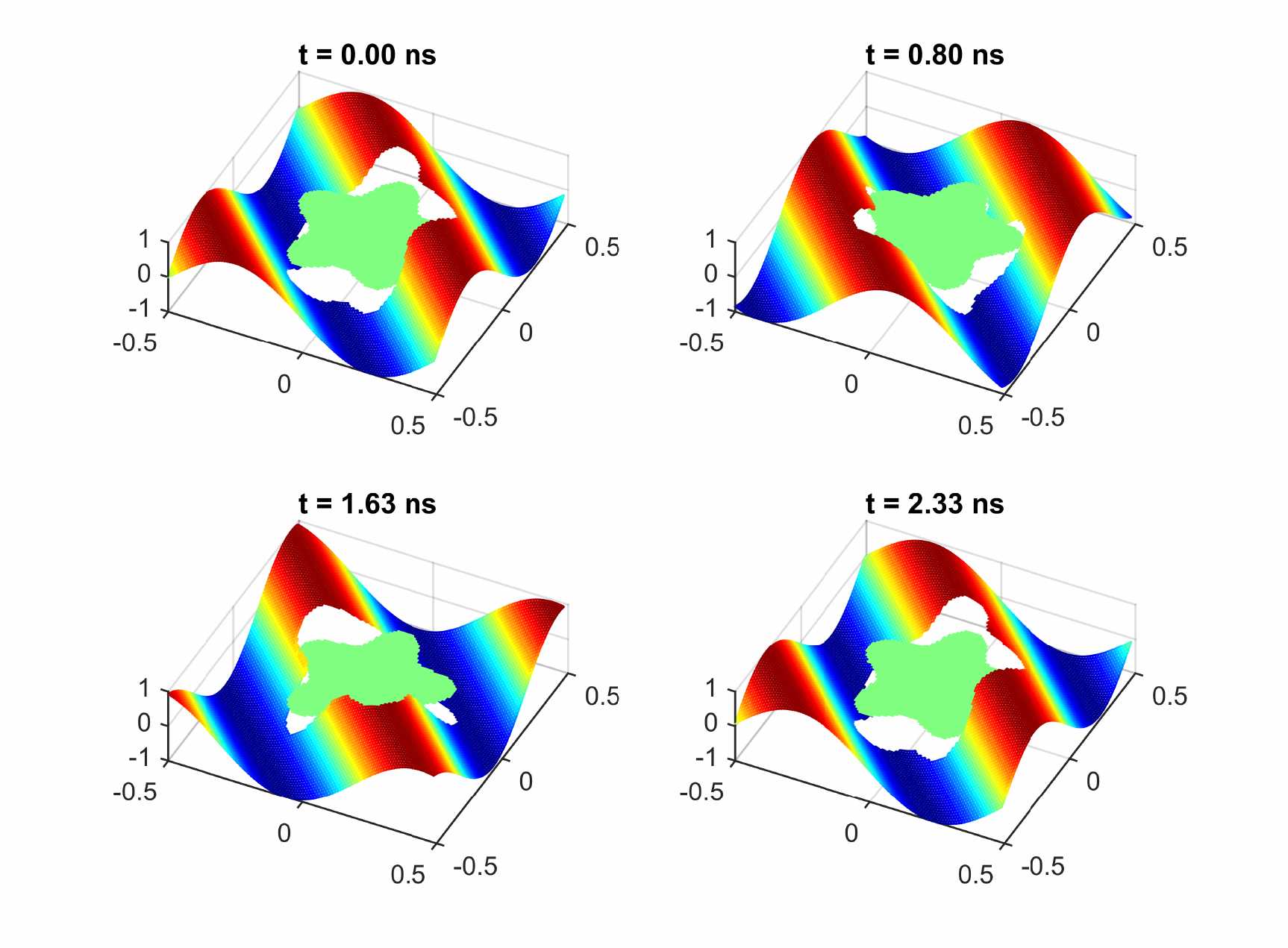}
\caption{Radiating interface solution at various time steps.}
\label{fig17}
\end{figure}

The solution was computed over one temporal period (approximately 2.359~ns), with a temporal discretization of $\Delta t = \frac{3}{4}\frac{\Delta x}{c}$. The resulting electric potential for a few time intervals, with $\Delta x = 0.01$, is plotted in figure \ref{fig17}. 

\begin{figure}[p]
\centering
\includegraphics[width=0.9\textwidth]{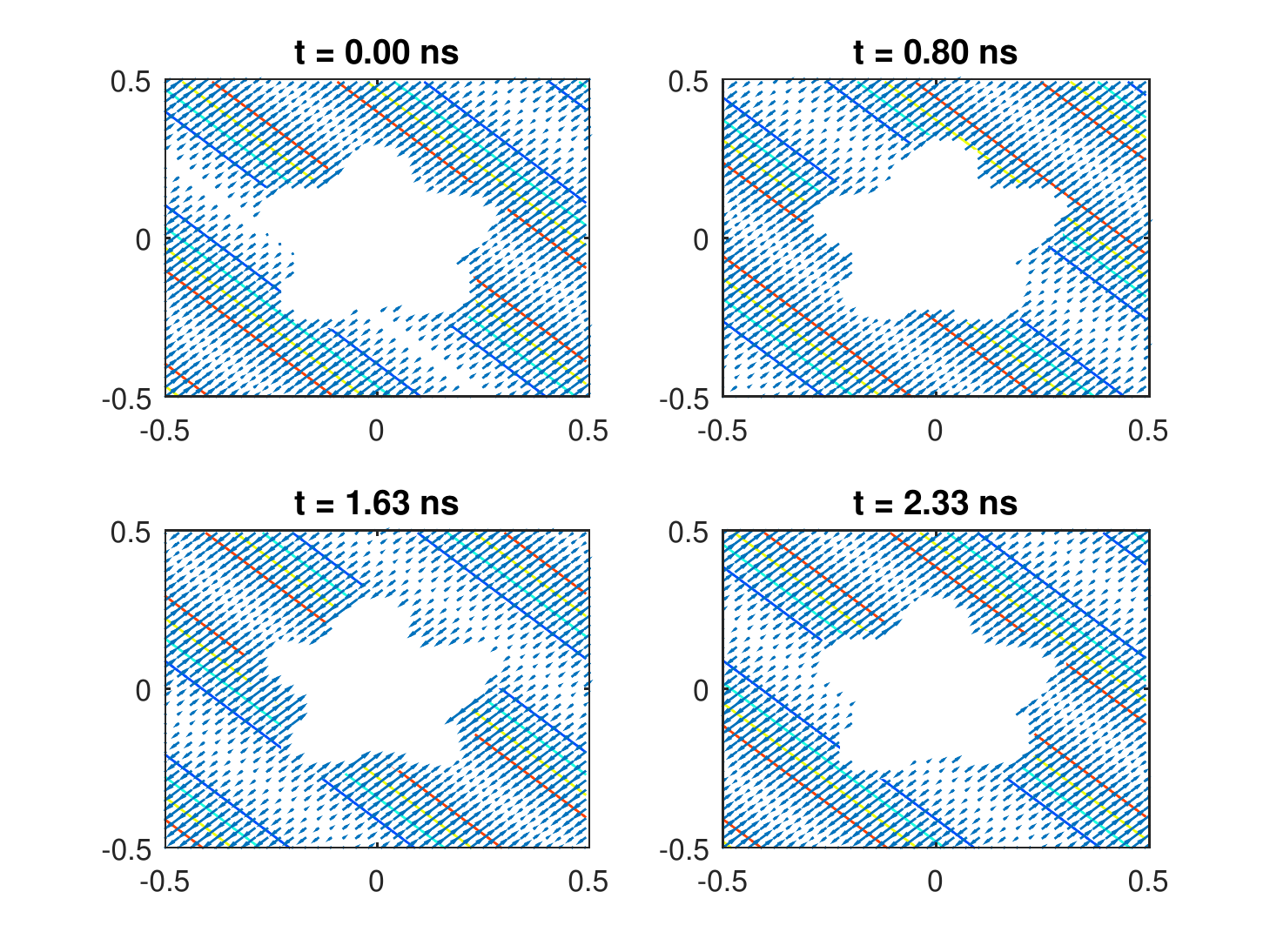}
\caption{Radiating interface solution gradient and equipotentials at various time steps.}
\label{fig18}
\end{figure}

\begin{figure}[t]
\centering
\includegraphics[width=0.9\textwidth]{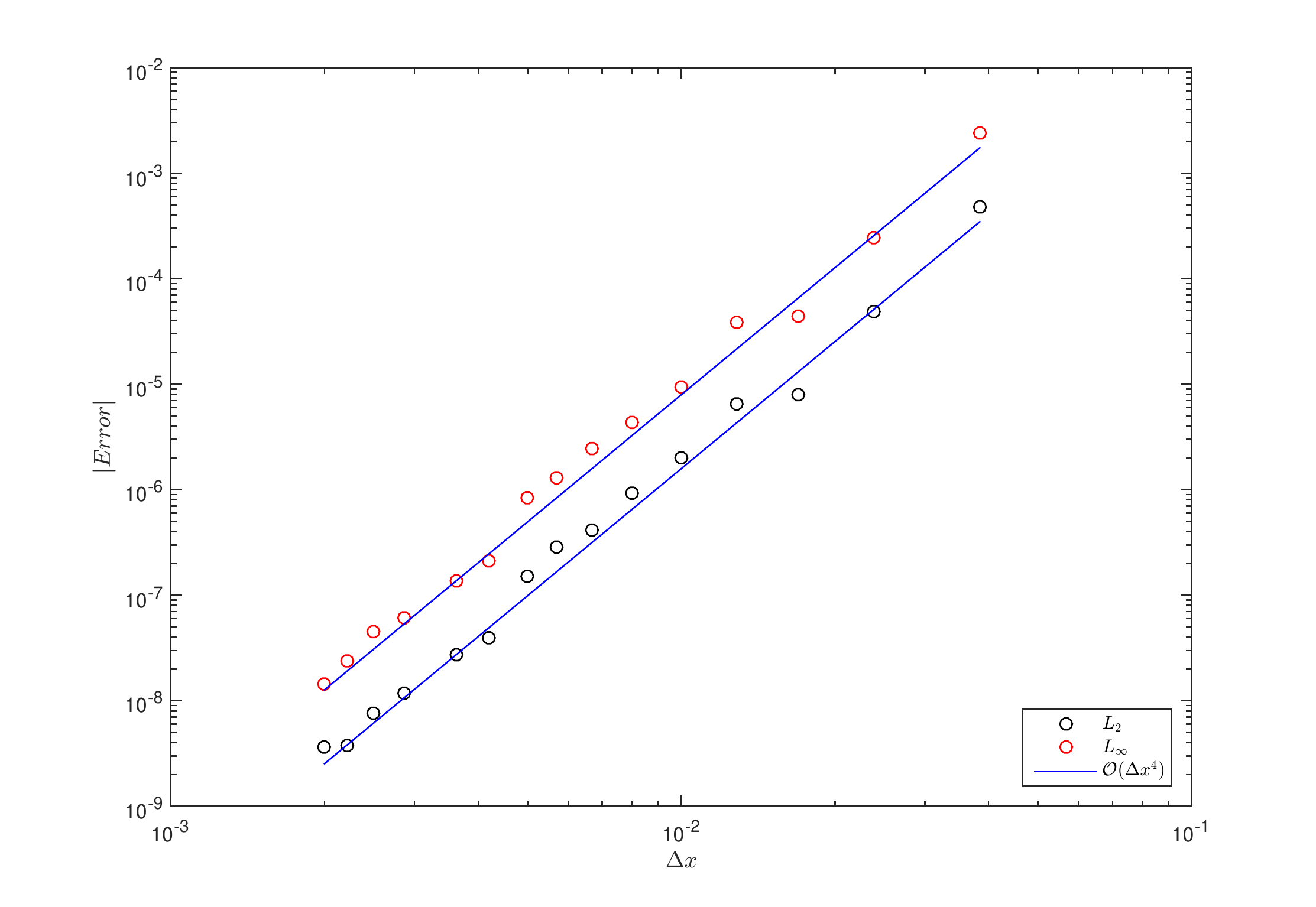}
\caption{2D Electromagnetic Scalar Potential - Solution error in $L_2$ and $L_\infty$ norms.}
\label{fig19}
\end{figure}

Moreover, the negative gradient of the potential, $-\nabla\phi$, representing the electrostatic contribution to the total $\vec{E}$ field of equation \eqref{E}, as well as equipotential contours of $\phi$, are plotted in figure \ref{fig18}. Lastly, the error as measured in the $L_2$ and $L_\infty$ norms is plotted in figure \ref{fig19} as a function of $\Delta x$, and once again exhibits the expected fourth order accuracy.

\subsection{Stability Verification}
With the method's accuracy thoroughly tested and demonstrated in the above implementations, data is now presented to validate the analysis of section \ref{stab}. The number of nodes along each spatial axis of the circular example in section \ref{circle}, and the 1D example in section \ref{line}, was varied between 100 and 500. In each case, the maximal value of the ratio $\gamma = \Delta t / \Delta x$ was measured for which the solution appeared to become unstable upon visual inspection (loss of smoothness, exponential increase, etc.). These values of $\gamma$ were then compared to the theoretically calculated values obtained through application of the standard stability region of RK4. The results are summarized in table \ref{tab1}.

\begin{table}[ht]
\centering
\begin{tabular}{|c|c|c|c|c|c|c|}
\cline{2-7}
\multicolumn{1}{r}{}& \multicolumn{3}{|c|}{1D} & \multicolumn{3}{|c|}{2D}\\
\hline
N & $\gamma_{t}$ & $\gamma_{c}$ & $\gamma_{CFM}$ & $\gamma_{t}$ & $\gamma_{c}$ & $\gamma_{CFM}$\\
\hline
100 & 1.23 & 1.31 & 1.24 & 1.23 & 1.31 & 1.26\\
\hline
200 & 1.23 & 1.26 & 1.24 & 1.23 & 1.26 & 1.24\\
\hline
300 & 1.23 & 1.25 & 1.23 & 1.23 & 1.25 &\multicolumn{1}{r}{}\\
\cline{1-6}
400 & 1.23 & 1.24 & 1.23 & 1.23 & 1.24 &\multicolumn{1}{r}{}\\
\cline{1-6}
500 & 1.23 & 1.24 & 1.23 &\multicolumn{3}{r}{}\\
\cline{1-4}
\end{tabular}
\caption{Comparison of theoretical and empirical stability limits for the CFM wave problem. $\gamma_t$ denotes the theoretically calculated maximum ratio of $\Delta t/\Delta x$ required for stability. $\gamma_c$ and $\gamma_{CFM}$ represent, respectively, the empirically determined stability limits for a continuous problem, and a problem with a jump discontinuity solved via the CFM.}
\label{tab1}
\end{table}

The $\gamma_{c}$ columns indicate the point of instability for a continuous problem in which no interface jumps occur and was determined by removing the circular interface, and allowing either $u^+$ or $u^-$ to occupy the whole domain. Moreover, $\gamma_{t}$ and $\gamma_{CFM}$ denote, respectively, the calculated theoretical and measured CFM instability points. From this data, it is clear that all empirically measured stability limits agree to good accuracy with the theoretically derived limit. Initial deviations can be seen in some of the data, however this is most likely attributed to the corresponding smaller values of $N$ (and therefore larger values of $\Delta x$) used at the start. As the grid is refined, the increased number of time steps required to span the same period aids in visually distinguishing any instability. Moreover, it also better approximates the asymptotic limit $\Delta x\rightarrow 0$, upon which the exact solution is based. These results thereby strongly corroborate the analysis presented in section \ref{stab}.

\section{Conclusion}

In this paper the Correction Function Method (CFM) has been used to obtain fourth order accurate solutions for the constant coefficient wave equation with interface jump conditions. While the present results were computed with fourth order accurate schemes and stencils, the method is, in principle, capable of being generalized to any arbitrarily high order. The CFM is similar to the Ghost Fluid Method (GFM), in the sense that it uses the concept of correction values applied at nodes affected by the discontinuity. Rather than assuming these corrections to be a discrete set, here the concept was extended to that of a correction function (CF), whose governing PDE and interface conditions were derived and defined within a small band surrounding the interface. Much like the GFM, the CFM retains the advantageous feature of being able to then model these corrections as equivalent source terms, allowing the use existing ``black box'' solvers with a simple augmentation to the system's right-hand-side.

In solving the correction function's governing PDE, a small band surrounding the interface was decomposed into square regions, in which the CF was expanded in terms of tricubic Hermite interpolants, and approximated via a least-squares minimization procedure. Since this process is independent of the relative position of the interface to the underlying grid, virtually any geometry is permitted, from interfaces with weak curvature grazing the grid, to those containing sharp cusps. The CF can thereby be fully precomputed and integrated into the selected non-compact five point spatial finite difference stencil as needed. RK4 was selected as the time marching method of choice, with necessary modifications to the CFM having been derived to maintain the solution's accuracy. These modifications represent a key novel contribution, as they demonstrate a systematic way to ensure accuracy and compatibility between, in principle, any time marching scheme and the CFM. The resulting fully fourth order scheme was applied to numerous problems in both one and two dimensions, resulting in all cases in clean, accurate, and sharp modeling of the interfaces, free of spurious oscillations, even in the vicinity of sharp corners.

The CFM is also an excellent candidate for parallelization, be it on Graphics Processing Units (GPUs) or other highly parallel architectures. The CF solution within a single space-time volume, spanning a node centered region and one time step, is fully independent and de-coupled from all other space-time volumes. As a result, each and every single region, for all time steps, may have its CF computed simultaneously with no communication between volumes. This allows maximum parallelizability, speed, and efficiency.

Lastly, it is strongly emphasized that the current method could, in principle, equally be used for problems in which the interfaces were moving, rather than static. Moreover, a great many interesting wave phenomena occur in domains for which the material parameters are not constant (such as the aforementioned dielectric interfaces in electromagnetics). While the method presented above is not directly applicable to such cases, a generalization of the wave CFM to these problems, similar to that proposed for the Poisson equation in \cite{[28]}, is currently under investigation. Furthermore, while the examples presented in section \ref{sec4} occurred in one or two spatial dimensions, in principle the provided method and procedure can be generalized to three dimensional problems.

\bibliographystyle{ieeetr}
\bibliography{References}

\begin{thebibliography}{10}

\bibitem{[1]}
A.~N. Marques, J.-C. Nave, and R.~R. Rosales, ``A correction function method
  for poisson problems with interface jump conditions,'' {\em Journal of
  Computational Physics}, vol.~230, pp.~7567--7597, Aug 2011.

\bibitem{[14]}
A.~N. Marques, {\em A Correction Function Method to Solve Incompressible Fluid
  Flows to High Accuracy With Immersed Geometries}.
\newblock PhD thesis, Massachusetts Institute of Technology, 2012.

\bibitem{[28]}
A.~N. Marques, J.-C. Nave, and R.~R. Rosales, ``High order solution of poisson
  problems with piecewise constant coefficients and interface jumps.''
\newblock arXiv preprint
  \href{https://arxiv.org/abs/1401.8084v3}{arXiv:1404.8084v3}, May 2016.

\bibitem{[6]}
C.~S. Peskin, ``Numerical analysis of blood flow in the heart,'' {\em Journal
  of Computational Physics}, vol.~25, pp.~220--252, Feb 1977.

\bibitem{[27]}
C.~Tu and C.~S. Peskin, ``Stability and instability in the computation of flows
  with moving immersed boundaries: A comparison of three methods,'' {\em SIAM
  Journal on Scientific and Statistical Computing}, vol.~13, pp.~1361--1376,
  Nov 1992.

\bibitem{[22]}
R.~J. LeVeque and Z.~Li, ``The immersed interface method for elliptic equations
  with discontinuous coefficients and singular source terms,'' {\em SIAM
  Journal on Numerical Analysis}, vol.~31, pp.~1019--1044, Aug 1994.

\bibitem{[25]}
Z.~Li, ``A fast iterative algorithm for elliptic interface problems,'' {\em
  SIAM Journal on Numerical Analysis}, vol.~35, pp.~230--254, Feb 1998.

\bibitem{[26]}
A.~Wiegmann and K.~P. Bube, ``The explicit-jump immersed interface method:
  Finite difference methods for pdes with piecewise smooth solutions,'' {\em
  SIAM Journal on Numerical Analysis}, vol.~37, no.~3, pp.~827--862, 2000.

\bibitem{[3]}
C.~Zhang, {\em Immersed Interface Methods for Hyperbolic Systems of Partial
  Differential Equations with Discontinuous Coefficients}.
\newblock PhD thesis, University of Washington, 1996.

\bibitem{[4]}
C.~Zhang and R.~J. LeVeque, ``The immersed interface method for acoustic wave
  equations with discontinuous coefficients,'' {\em Wave Motion}, vol.~25,
  pp.~237--263, May 1997.

\bibitem{[2]}
S.~Deng, ``On the immersed interface method for solving time-domain maxwell's
  equations in materials with curved dielectric interfaces,'' {\em Computer
  Physics Communications}, vol.~179, pp.~791--800, Dec 2008.

\bibitem{[19]}
B.~Lombard and J.~Piraux, ``Numerical treatment of two-dimensional interfaces
  for acoustic and elastic waves,'' {\em Journal of Computational Physics},
  vol.~195, pp.~90--116, Mar 2004.

\bibitem{[20]}
R.~P. Fedkiw, T.~Aslam, B.~Merriman, and S.~Osher, ``A non-oscillatory eulerian
  approach to interfaces in multimaterial flows (the ghost fluid method),''
  {\em Journal of Computational Physics}, vol.~152, pp.~457--492, Jul 1999.

\bibitem{[15]}
J.~Piraux and B.~Lombard, ``A new interface method for hyperbolic problems with
  discontinuous coefficients: One-dimensional acoustic example,'' {\em Journal
  of Computational Physics}, vol.~168, pp.~227--248, Mar 2001.

\bibitem{[23]}
J.~D. Jackson, {\em Classical Electrodynamics}.
\newblock Hoboken, New Jersey: John Wiley \& Sons, Inc., third~ed., 1998.

\bibitem{[24]}
D.~J. Griffiths, {\em Introduction to Electrodynamics}.
\newblock Upper Saddle River, New Jersey: Prentice Hall, third~ed., 1999.

\bibitem{[8]}
N.~K. Georgieva and H.~W. Tam, ``Potential formalisms in electromagnetic-field
  analysis,'' {\em IEEE Transactions on Microwave Theory and Techniques},
  vol.~51, pp.~1330--1338, Apr 2003.

\bibitem{[9]}
N.~Georgieva and E.~Yamashita, ``Time-domain vector-potential analysis of
  transmission-line problems,'' {\em IEEE Transactions on Microwave Theory and
  Techniques}, vol.~46, pp.~404--410, Apr 1998.

\bibitem{[13]}
R.~Dyczij-Edlinger and O.~Biro, ``A joint vector and scalar potential
  formulation for driven high frequency problems using hybrid edge and nodal
  finite elements,'' {\em IEEE Transactions on Microwave Theory and
  Techniques}, vol.~44, pp.~15--23, Jan 1996.

\bibitem{[17]}
O.~Biro, ``Edge element formulations of eddy current problems,'' {\em Computer
  Methods in Applied Mechanics and Engineering}, vol.~169, pp.~391--405, Feb
  1999.

\end{thebibliography}
\end{document}